\documentclass[12pt]{article}
\usepackage[utf8]{inputenc} 
\usepackage[T1]{fontenc} 
\usepackage{amsmath}
\usepackage{amsthm}
\usepackage{amsfonts}
\usepackage{amssymb}
\usepackage{graphicx} 
\usepackage{latexsym}
\usepackage{color}

\advance \topmargin by -\headheight
\advance \topmargin by -\headsep
\evensidemargin \oddsidemargin
\newtheorem{theorem}{Theorem}[section]

\newtheorem{lemma}[theorem]{Lemma}

\theoremstyle{definition}
\newtheorem{definition}[theorem]{Definition}

\newtheorem{example}[theorem]{Example}

\newtheorem{claim}[theorem]{Claim}
\theoremstyle{remark}
\newtheorem{remark}[theorem]{Remark}

\numberwithin{equation}{section}

\newcommand{\GL}{\mathrm{GL}}

\newcommand{\SL}{\mathrm{SL}}

\newcommand{\SU}{\mathrm{SU}}

\newcommand\Lie[1]{{\mathfrak{#1}}}

\newcommand\CC{\mathbb C}

\newcommand\HH{\mathbb H}
\newcommand\NN{\mathbb N}
\newcommand\RR{\mathbb R}
\newcommand\ZZ{\mathbb Z}
\newcommand\QQ{\mathbb Q}

\newcommand\KK{\mathbb K}

\newcommand\FF{\mathbb F}

\newcommand\cA{\mathcal{A}}
\newcommand\cB{\mathcal{B}}

\newcommand\cD{\mathcal{D}}

\newcommand\cH{\mathcal{H}}

\newcommand\cW{\mathcal{W}}
\newcommand\cE{\mathcal{E}}
\newcommand\cU{\mathcal{U}}

\newcommand\cK{\mathcal{K}}

\newcommand\cJ{\mathcal{J}}

\newcommand{\ma}{\mathfrak{a}}

\newcommand\norm[1]{\left\|#1\right\|}

\newcommand\abs[1]{\left|#1\right|}

\newcommand\inn[1]{\left\langle #1 \right\rangle}
\newcommand\set[1]{\left\{{#1}\right\}}

\DeclareMathOperator{\vol}{vol}

\DeclareMathOperator{\Hom}{Hom}

\DeclareMathOperator{\re}{Re}

\DeclareMathOperator{\supp}{supp}
\DeclareMathOperator{\mm}{m}

\begin{document}
\title{Automorphic density estimates and optimal approximation exponents}

\author{Mikołaj Frączyk, Alexander Gorodnik and Amos Nevo} 
%
%
%

%
%
%




\maketitle
\begin{abstract}
The present paper is devoted to establishing an optimal  approximation exponent for the action of an irreducible uniform lattice subgroup of a product group on its proper factors. 
Previously optimal approximation exponents for lattice actions on homogeneous spaces were established under the assumption that the restriction of the automorphic representation to the stability group is suitably tempered. 
However, for irreducible lattices in semisimple algebraic groups, either this property does not hold or it amounts to an instance of the 
Ramanujan-Petersson-Selberg conjecture.
Sarnak's Density Hypothesis and its variants bounding the multiplicities of irreducible representations occurring in the decomposition of the automorphic representation can be viewed as a weakening of the temperedness property. A refined form of this hypothesis has recently been established for  uniform irreducible arithmetic congruence lattices arising from quaternion algebras. We employ this result in order to establish - unconditionally - an optimal approximation exponent for the actions of these lattices on the associated symmetric spaces. We also give a general spectral criterion for the optimality of the approximation exponent for  irreducible uniform lattices in a product of arbitrary Gelfand pairs. Our methods involve utilizing the multiplicity bounds in the pre-trace formula, establishibng refined estimates of the spherical transforms, and carrying out an elaborate spectral analysis that bounds the Hilbert--Schmidt norms of carefully balanced geometric convolution operators. 
\end{abstract}

\tableofcontents

\section{Introduction }

\subsection{The approximation exponent problem for  homogeneous spaces}
This paper is devoted to some important instances of the following fundamental problem regarding quantitative density of orbits on homogeneous spaces.
Let $G$ be a locally compact second countable non-compact group, $X$ a homogeneous space of $G$, and $\Gamma$ a discrete subgroup of $G$.
We consider the action of $\Gamma$ on $X$ (which we assume to have  at least one dense orbit) and aim to establish quantitative gauges measuring the density of the corresponding $\Gamma$-orbits. For this purpose, we fix a metric $D$ on $X$ and a proper function $\abs{\,\cdot\,} :G\to \RR^+$.
Given a dense orbit $\Gamma x\subset X$ and a point $y\in X$, we will be interested in existence of solutions $\gamma\in \Gamma$ for the inequalities
\begin{equation}\label{eq:ddd}
D(\gamma\, x,y)\le \epsilon \quad \hbox{and}\quad \abs{\gamma}\le R,
\end{equation}
where $\epsilon\to 0$ and $R\to\infty$ with prescribed rates.
Let us suppose that there exist $\zeta>0$ and $\epsilon_0=\epsilon_0(x,y,\zeta)>0$ such that for every $\epsilon\in (0,\epsilon_0)$,
the system of inequalities 
$$
D(\gamma \, x,y)\le \epsilon \quad \hbox{and}\quad \abs{\gamma}\le \zeta\, \log(1/\epsilon).
$$
has solutions $\gamma\in \Gamma$. Following \cite[Def.~2.1]{GGN14}, we define the \emph{approximation exponent} $\kappa_\Gamma(x,y)$
as the infimum of $\zeta$ such that the foregoing property holds.
The exponent $\kappa_\Gamma(x,y)$ constitutes a natural gauge of the efficiency of approximating the point $y$ by elements in the orbit $\Gamma  x$.  

An extensive class of Diophantine approximation problems, and in particular the approximation exponent, were originally raised by Lang \cite{L65}, in the setting of actions 
of arithmetic groups on homogeneous varieties of linear algebraic groups.  
These problems include analogs of the classical Dirichlet's, Khinchine's and Schmidt's theorems in  Diophantine approximation, extended to the context of lattice orbits on homogeneous spaces. A systematic approach to investigate these problem has been developed in recent years, and for a detailed discussion of it we refer to \cite{GGN12, GGN13, GGN14, GN23a}. In the present paper, however, our focus will be on establishing optimal results for the approximation exponent. 

Let us assume further that $G$ is a unimodular group and  $X=G/H$ for a closed unimodular subgroup $H$. 
We denote the invariant measures on the groups $G$ and $H$ by $m_G$ and $m_H$, respectively, and the $G$-invariant measure on the space $X$ by $m_X$.
There is a natural upper bound governing  how many approximations as in \eqref{eq:ddd} can be constructed for a generic point $y$. Namely, the number of solutions of \eqref{eq:ddd} can be expected to be bounded by the number of 
points of the orbit $\Gamma  x$ which are contained in compact neighbourhoods of $y$ in $X$. Let us assume, for simplicity, that the limit
$$
a(H):=\lim_{R\to\infty} \frac{1}{R}\log m_H\big(\{h\in H:\, \abs{h}\le R\}\big)
$$
exists and is positive. Under certain mild regularity conditions on the gauge function $\abs{\,\cdot\,}$ and the metric $D$,
one can show that the number of points $\gamma \, x$ lying in a compact subset $\Omega\subset X$ with $\gamma\in \Gamma$ satisfying $\abs{\gamma}\le R$ is 
at most $O_{\Omega,x,\delta} \big(e^{(a(H)+\delta)R}\big)$, with $\delta>0$ arbitrarily small.
This upper bound on the number of points available for the approximation can be expected to constrain the quality of approximations by elements of $\gamma\in \Gamma$ satisfying $\abs{\gamma}\le R$.
Developing this idea, it was shown in \cite[Th.~3.1]{GGN14} that under the assumption of ergodicity for the action of  $\Gamma$ on $X$,  for every $x\in X$ with a dense $\Gamma$-orbit and for almost every $y\in X$, the following pigeonhole lower bound holds: 
\begin{equation}\label{eq:low}
\kappa_\Gamma(x,y)\ge \frac{\dim(X)}{a(H)}.
\end{equation}
Therefore, the challenge is to establish a comparable upper bound for $\kappa_\Gamma(x,y)$.

\medskip

When $\Gamma$ is a lattice subgroup in $G$, an approach to establishing upper bounds for the approximation exponents $\kappa_\Gamma(x,y)$
was developed in \cite{GGN12,GGN14} using spectral estimates in the automorphic representation $\lambda_{\Gamma \backslash G}$
of $G$ on $L^2(\Gamma \backslash G)$.
The key to the method is to establish a {quantitative mean ergodic theorem} for the action of the subgroup $H$ on the space $\Gamma\backslash G$
equipped with the invariant probability measure $m_{\Gamma\setminus G}$ induced by the measure $m_G$. 
This provides an effective speed of distribution for $H$-orbits in $\Gamma\backslash G$, and allows the use of dynamical considerations in the form of a "shrinking target argument".  These considerations are combined with a quantitative duality argument,  which translates an equidistribution rate for $H$-orbits on $\Gamma\backslash G$ to a quantitative density estimates for the $\Gamma$-orbits in $X=G/H$. 

To formulate the relevant spectral estimate, we set
$$
H_R:=\set{h\in H\,:\,  |h|\le  R}
$$
and introduce the corresponding averaging operators
$$
\lambda_{\Gamma\backslash G}(\beta_R) :L^2(\Gamma\setminus G)\to L^2(\Gamma\setminus G)
$$
defined by
$$
\lambda_{\Gamma\backslash G}(\beta_R)\phi(z):=\frac{1}{m_H(H_R)}\int_{H_R}\phi(zh)\,dm_H(h)\quad \hbox{for $\phi\in L^2(\Gamma\setminus G)$ and $z\in \Gamma\backslash G\,.$}
$$
The \emph{quantitative mean ergodic theorem} for the operators $\lambda_{\Gamma\backslash G}(\beta_R)$ states that  there exists $\theta > 0$ such that 
\begin{equation}\label{eq:mean}
\left\|\lambda_{\Gamma\backslash G}(\beta_R)\phi-\int_{{\Gamma\backslash G}} \phi\, dm_{\Gamma\backslash G} \right\|_{L^2(\Gamma\setminus G)}\ll_\delta m_H(H_R)^{-\theta+\delta}\,\|\phi\|_{L^2(\Gamma\backslash G)},\quad \phi\in L^2(\Gamma\backslash G),
\end{equation}
for every $\delta > 0$ and sufficiently large $R$. Let us denote by $\theta=\theta(G,H,\Gamma)$ the maximal $\theta$ for which  
property \eqref{eq:mean} holds. 
It was shown in \cite[Th.~3.3]{GGN14} that for almost every $x\in X$ and every $y\in X$,
\begin{equation} \label{eq:up}
\kappa_\Gamma(x,y)\le \frac{1}{2 \theta(G,H,\Gamma)}\cdot \frac{\dim(X)}{a(H)}.
\end{equation}
In particular, when $\theta(G,H,\Gamma)=1/2$, the lower bound \eqref{eq:low} coincides with the upper bound \eqref{eq:up}, establishing the optimal approximation exponents for almost every pair $(x,y)\in X^2$.

We note that the estimate \eqref{eq:mean} with $\theta=1/2$ holds 
when  $H$ is a simple non-compact Lie group, and the representation of $H$ on the space $L^2_0(\Gamma\backslash G)$, the space of functions with zero integral,
is weakly contained in the regular representation of $H$. 
In this case the required estimate \eqref{eq:mean} with $\theta=1/2$ is a reflection of the \emph{Kunze-Stein phenomenon} \cite{c78}.
As elaborated in \cite{GGN14}, the condition $\theta(G,H,\Gamma)=1/2$ holds for an extensive class of triples $(G,H,\Gamma)$. 
However, it is also known that $\theta(G,H,\Gamma)>1/2$ 
in many cases, and in some cases
the temperedness condition $\theta(G,H,\Gamma)=1/2$ amounts to a version of the Ramanujan-Petersson-Selberg conjecture,  so that the temperedness condition conjecturally holds, but is currently out of reach.

This brings to the fore the natural problem of whether the 
estimate \eqref{eq:up} can be improved. In \cite{GN23b}, 
the case when the group $\Gamma$ is an $S$-arithmetic lattice in simple algebraic group projecting densely on the Archimedean factor $G_\infty$ was considered, and it was shown (see \cite[Theorem 1.4]{GN23b}) that 
estimates on the $L^2$-discrepancy describing the distribution of $\Gamma$-points in $G_\infty$ imply norm bounds 
for the corresponding automorphic representations.
In particular, the discrepancy estimate corresponding to 
the pigeonhole bound implies that the discrete part of the spectrum 
is tempered.  Remarkably, despite this result, 
it turns out that the optimal bound on the approximation exponent 
can be established in some cases without knowing that 
$\theta(G,H,\Gamma)=1/2$.
Sardari \cite{S19} established the optimal density
for rational points on quadratic surfaces in at least five
variables using a refined analysis based on the Circle Method.
Jana and Kamber \cite{JK22} investigated 
quantitative density  of orbits of 
$\Gamma=\hbox{SL}_d(\ZZ[1/p])$
acting on the symmetric space 
$X=\hbox{SL}_d(\RR)/\hbox{SO}_d(\RR)$.
They proved optimal bounds on 
the approximation exponent $\kappa_\Gamma(x,y)$, for generic $(x,y)\in X^2$, for $d=2,3$ and for $d\ge 4$ a conditional optimal result, assuming the Sarnak Density Conjecture
 discussed in Section \ref{sec:density} below.
The key idea of the proof in \cite{JK22} is to exploit
bounds on the multiplicities of non-tempered irreducible representations appearing in the spectral decomposition of $L^2\big(\Gamma\backslash (\hbox{SL}_d(\RR)\times \hbox{SL}_d(\QQ_p)) \big)$. 

Let us now discuss a problem described by Sarnak as the "non-Archimedian" counterpart of  the optimal diophantine exponent, namely the optimal lifting problem for finite quotients, where the idea of using multiplicity bounds to establish optimality originated. 

Let $\Gamma$ be a discrete group equipped with a proper function $|\cdot|:\Gamma\to \mathbb{R}^+$ and $\Gamma_n$  a family of finite index subgroups of $\Gamma$ such that $|\Gamma/\Gamma_n|\to \infty$. One considers the factor maps $\pi_n: \Gamma\to \Gamma/\Gamma_n$.
For $x\in \Gamma/\Gamma_n$, one seeks solutions $\gamma\in\Gamma$
satisfying 
\begin{equation}\label{eq:lift}
\pi_n(\gamma)=x\quad \hbox{and}\quad |\gamma|\le R
\end{equation}
with $n \to\infty$ and $R\to \infty$.
Sarnak \cite{S15} considered this problem for $\Gamma=\hbox{SL}_2(\ZZ)$
and $\Gamma_n$'s being the congruence subgroups.
He showed that generically the lifting problem
has solutions satisfying $\|\gamma\|_\infty\le |\Gamma/\Gamma_n|^{1/2+\epsilon}$ for all $\epsilon>0$ and sufficiently large $R$,
which is essentially optimal.
While this result can be deduced from the Ramanujan--Petersson--Selberg Conjecture, which can be 
equivalently formulated as $\theta(\hbox{SL}_2(\RR),\hbox{SL}_2(\RR),\Gamma_n)=1/2$,
remarkably  Sarnak's result shows that this conjecture can be circumvented in 
establishing the optimal lifting property.
Subsequently, Golubev and Kamber \cite{GK23} developed a general framework for studying the above lifting problem, which
allows to deduce existence of solutions for \eqref{eq:lift}
with prescribed relation between $|\Gamma/\Gamma_n|$
and $R$ from bounds on multiplicities for non-tempered 
irreducible representations appearing in the spectral
decompositions of $L^2(G/\Gamma_n)$. 
In particular, they proved that generically the optimal lifting property can be deduced from a version of the Sarnak Density
Conjecture.
Recently, the combined work of Jana and Kamber \cite{JK24} and Assing and Blomer \cite{AB23} established a suitable version of the Sarnak Density Conjecture for the principal congruence subgroups  of $\hbox{SL}_d(\ZZ)$ which suffices to prove the optimal lifting property in this setting. We also refer to \cite{KL23} where the optimal 
lifting property has been established for another class of 
finite-index subgroups of $\hbox{SL}_3(\ZZ)$.

Our two main goals in the present paper are as follows. 
First, we give a spectral and geometric criterion 
(Theorem \ref{th:general} below) which suffices to establish optimality of the approximation exponent in a very general set-up when the group $\Gamma$ is an irreducible uniform lattice in the product $G_1\times G_2$ of
locally compact groups. 
Second, we utilize the recently established generalization of the  Sarnak Density Conjecture in the context of uniform arithmetic congruence lattices arising from general quaternion algebras \cite{FHMM} and derive the optimal approximation exponent (Theorem \ref{th:quaternion} below) for their action on the associated symmetric spaces. Namely, we prove unconditional spectral estimates serving as a substitute for the temperedness condition, sidestepping any appeal to the Ramanujan--Petersson--Selberg conjecture.  

As an example, our analysis applies when $\Gamma$ is a uniform irreducible arithmetic congruence lattice in $G=\hbox{SL}_2(\RR)\times 
\hbox{SL}_2(\RR)$, arising from a quaternion algebra, with $H=\hbox{SO}_2(\RR)\times \hbox{SL}_2(\RR)$, and $G/H\simeq\HH^2$. 
Note that temperedness of the representation of $H$ on $L_0^2(\Gamma\backslash G)$ (which is sufficient to derive the optimal exponent for $\Gamma$ in its action on $X=G/H$) certainly requires a version of the 
Ramanujan-Petersson-Selberg conjecture. However, we will prove the optimality of the Diophantine exponent unconditionally, using 
the multiplicity estimates established in \cite{FHMM}. 
We note that important progress towards the Sarnak density conjecture for uniform arithmetic congruence lattices in $PSL_2(\RR)^d$ was established previously by Kelmer \cite{Ke11}. While the uniform estimates our method employs are not stated  explicitly in \cite{Ke11}, it is possible that they can be derived from the arguments there.

\subsection{Method of proof in the present paper}

\subsubsection{Trace estimates of geometric convolution operators} 

We now describe briefly the approach we take to estimating the approximation exponents, which is applicable in the following general setting. Let 
$G=G_1\times G_2$ be a product group, and $K_i\subset G_i$
are compact subgroups such that $(G_i,K_i)$ are Gelfand pairs for $i\in \set{1,2}$.
We set $X_i:=G_i/K_i$.
Let $\Gamma$ be a uniform irreducible lattice in $G$.
We consider the action of $\Gamma$ on $X_1$, via its dense projection to the first factor. Let us fix proper $G_i$-invariant metrics $D_i$ on $G_i/K_i$ for $i\in \set{1,2}$.
For $x_1,y_1\in X_1$,
we consider the inequalities
\begin{equation}\label{eq:ineqq}
D_1(\gamma_1 x_1,y_1)\le r\quad \hbox{and} \quad D_2(\gamma_2 K_2,K_2)\le R
\end{equation}
with $\gamma=(\gamma_1,\gamma_2)\in \Gamma$.
In this case, the volume growth rate is given by
$$
 a(G_2):=\lim_{R\to\infty} \frac{1}{R}\log m_{G_2}\big(\{g_2\in G_2:\, D_2(g_2 K_2,K_2)\le R\}\big).
$$
The main ingredient in our proof of an upper bound for the approximation exponent will rely on spectral estimates of  suitable geometrically defined operators acting in the automorphic representation $L^2(\Gamma \backslash G)$. However, in the present approach the operators will have supports with non-empty interiors in $G$ (rather than be supported in $H$, as in the previous approach indicated above). They will then be trace-class operators and we will estimate their Hilbert-Schmidt norms using bounds on their traces.  
In the special case $G$ is a connected Lie group, our main result (Theorem \ref{th:general}) states that if a suitable trace bound holds,  then 
for all $\zeta> \frac{\dim(X_1)}{a(G_2)}$ and 
almost every $x_1,y_1\in X_1$, the inequalities 
$$
D_1(\gamma_1 x_1,y_1)\le r\quad \hbox{and} \quad D_2(\gamma_2 K_2,K_2)\le \zeta \, \log(1/r)
$$
have solutions $\gamma=(\gamma_1,\gamma_2)\in \Gamma$ for all sufficiently small $\epsilon>0$. Therefore, 
\begin{equation}\label{eq:oopt}
\kappa_\Gamma(x_1,y_1)\le \frac{\dim(X_1)}{a(G_2)}
\quad\hbox{for almost every $(x_1,y_1)\in X_1\times X_1$,}
\end{equation}
which matches precisely with the lower bound \eqref{eq:low},
so that this gives the exact value for the exponent $\kappa_\Gamma(x_1,y_1)$ generically, i.e. for almost every pair.

To give an idea of the spectral trace bound that we employ, let us introduce functions 
$$
F_{r,R}:=f_{1,r}\otimes f_{2,R}\quad\hbox{with 
$f_{1,r}\in C_c(G_1)$ and $f_{2,R}\in C_c(G_2)$}
$$
satisfying
\begin{align*}
\supp(f_{1,r}) &\subset \widetilde{B}_r^{G_1}:=\{g\in G_1:\, D_1(gK_1,K_1)\le r\},\\
\supp(f_{2,R}) &\subset \widetilde{B}_R^{G_2}:=\{g\in G_2:\, D_2(gK_2,K_2)\le R\}.
\end{align*}
We consider the corresponding convolution operators defined by
\begin{equation}\label{eq:av0}
\lambda_{\Gamma\backslash G}(F_{r,R}):L^2(\Gamma\backslash G)\to L^2(\Gamma\backslash G): \phi\mapsto \int_{G}\phi(zg)F_{r,R}(g)\,dm_{G}(g)
\end{equation}
and establish a bound on the Hilbert-Schmidt norm of $\lambda_{\Gamma\backslash G}(F_{r,R})$.
Analyzing $\lambda_{\Gamma\backslash G}(F_{r,R})\phi$
for suitably chosen test functions $\phi \in L^2(\Gamma\backslash G)$, we derive existence of solutions for \eqref{eq:ineqq}. 
It turns out that the optimal approximation exponent
corresponds to the choice of parameters $r$ and $R$
according to the {\it matching volumes condition}:
$$
m_{G_1}\big(\widetilde{B}_r^{G_1}\big)\cdot m_{G_2}\big(\widetilde{B}_R^{G_2}\big)\asymp 1.
$$

We note that 
while the approach of \cite{GGN14} was based on the operators defined as in \eqref{eq:mean} and depended on the automorphic spherical representation restricted to $H$ being tempered, the present approach affords dealing with the part of the spectral decomposition of the representation $\lambda_{\Gamma\backslash G}$ where the norms of the convolution operators
have slower decay. This is achieved by surmounting the challenge of explicitly bounding  the trace in (\ref{eq:av0}), which relies on estimates bounding the multiplicity of the appearance of representations with slow decay, namely on density estimate in the automorphic representation.  

We execute this approach for uniform congruence
arithmetic lattices arising from quaternion algebras
(see Section \ref{sec:quaternions} for the definitions).
For instance, if $\Gamma\subset \hbox{SL}_2(\RR)\times \hbox{SL}_2(\RR)$ is a uniform 
congruence arithmetic lattice in $\hbox{SL}_2(\RR)\times \hbox{SL}_2(\RR)$, our main result (Theorem \ref{th:quaternion}) implies that the optimal exponent is $\kappa(x,y)=1$, namely for  
any $\zeta> 1$  the inequalities 
$$
d(\gamma_1 \, x,y)\le \epsilon\quad \hbox{and} \quad d(\gamma_2 \, i,i)\le \zeta \, \log(1/\epsilon)
$$
have solutions $\gamma=(\gamma_1,\gamma_2)\in \Gamma$
for almost every $x,y\in \HH^2$ and sufficiently small $\epsilon$.

\subsubsection{Density estimates}\label{sec:density}

The original Density Conjecture is due to Sarnak and Xue and appeared in \cite{SX91}.
Let us give an account of a general density estimate for a semi-simple real Lie group $G$. Let $\mathcal H$ be a Hilbert space with an irreducible unitary representation $\pi\colon G\to \mathcal U(\mathcal H)$. A matrix coefficient of $\pi$ is the function $g\mapsto \langle v, \pi(g)w\rangle$, with $v,w\in \mathcal H$. We define a numerical invariant $p(\pi)$, measuring the extent to which $\pi$ fails to be tempered:
$$p(\pi):=\inf\{p\geq 2\mid \pi \text{ has a non-zero matrix coefficient in } L^{p}(G)\}.$$

A representation $\pi$ is tempered if and only if $p(\pi)=2$. The original formulation of the density conjecture pertained to the level aspect, asserting that for any arithmetic lattice $\Gamma\subset G$, any sequence of principal congruence subgroups $\Gamma_n\subset \Gamma$ and any irreducible unitary representation $\pi$ and $\varepsilon>0$, we have 
$${\rm m}(\pi,\Gamma_n):=\dim \Hom(\pi, L^2(\Gamma_n\backslash G))\ll_{\varepsilon} \vol(\Gamma_n\backslash G)^{2/ p(\pi)+\varepsilon}.$$

The obvious bound in this context would be  linear in the volume $\vol(\Gamma\backslash G)$. Therefore, the Density Hypothesis predicts a power  saving whose strength should only depend on the invariant $p(\pi).$ The conditions for the sequence of lattices can be relaxed considerably, for example by dropping the assumption that all the lattices $\Gamma_n$ are contained a single ambient lattice (see \cite{FHMM} for an instance of density estimate in this generality).
In practice, we seek a bound that counts the total multiplicity of the irreducible representations belonging to a specific spectral window which might also vary (we refer to Section \ref{sec:densityestimate} for more details). The volume on the right-hand side of the density estimate should be replaced by a quantity that depends on the spectral window and represents a bound which does not depend on the non-tempered parameters. In this way, one can talk about non-trivial density estimate for a fixed lattice $\Gamma,$ where instead of varying the lattice group (i.e. the level aspect) we vary the archimedean parameters of the representation $\pi$. This archimedean aspect plays the key role in the present paper. 
For the state of the art results and more in depth overview of the Density Hypotheses we refer the reader to \cite{AB23, FHMM}.


\subsection{Organization of the paper}
In Section \ref{sec:statement} we setup basic notation
and state our main results (Theorem \ref{th:general} and Theorem \ref{th:quaternion}).
Theorem \ref{th:general} is proved in Section
\ref{sec:proof_general},
which involves relating the approximation problem
to norm bounds for suitable integral operators in Section
\ref{sec:geom} and analyzing the corresponding norms in Section 
\ref{sec:spec}.
In Section \ref{sec:quaternions} we
review the construction of arithmetic lattices
based on quaternion algebras.
The proof of Theorem \ref{th:quaternion} is completed in Section \ref{sec:proof_quat}.
In particular, in Section \ref{sec:spherical functions}
we construct a suitable family of
test-functions with controlled spherical transforms.
In Section \ref{sec:densityestimate}
we review basic properties of spectrum of the relevant homogeneous spaces, 
and in Section \ref{sec:trace} we establish a crucial trace estimate,
which allows us to reduce the proof of Theorem \ref{th:quaternion}
to Theorem \ref{th:general}.
Finally, in Section \ref{sec:local SF} we deal with local
estimates on spherical transforms (Theorem \ref{lem-1}),
which is needed in Section \ref{sec:spherical functions}.

\subsection*{Notation} 

We write $A\ll_{\epsilon,\delta,\ldots} B$ if there exists a constant $c>0$ depending only on the parameters $\epsilon,\delta,\ldots$ such that $A\le c\, B$
and $A\asymp_{\epsilon,\delta,\ldots}B$ if there exists constants $c_1,c_2>0$ depending only on the parameters $\epsilon,\delta,\ldots$ such that $c_1\, B\le A\le c_2\, B$.
 We will follow the  "variable constant convention" where the value of a constant  may be different in the various estimates appearing in the same argument.

\subsection*{Acknowledgements}
A.G. was supported by SNF Grant 200020--212617. M. F. acknowledges the support from the Dioscuri programme initiated by the Max Planck Society, jointly managed with the National Science Centre in Poland, and mutually funded by Polish the Ministry of Education and Science and the German Federal Ministry of Education and Research.
A. N. was partially supported by ISF Moked Grant 2019-19. 
\section{Preliminaries and statements of main results}
\label{sec:statement}

\subsection{The pre-trace formula and Gelfand pairs}
Let $G$ be a unimodular locally compact second countable non-compact group, and let $m_G$ denote a (left and right) Haar measure on $G$. We denote by $\widehat{G}$ the set of equivalence classes of irreducible unitary representations of $G$, equipped with the Fell topology.

Given any strongly continuous unitary representation $\pi : G \to \cU(\cH_\pi)$ to  the unitary group of a separable Hilbert space $\cH_\pi$, any function $F\in L^1 (G)$  defines a bounded operator $\pi(F) : \cH_\pi\to \cH_\pi$, given by
$$\inn{\pi(F)u,v}_{\cH_\pi}:=\int_G F(g)\inn{\pi(g)u,v}_{\cH_\pi} dm_G(g),\quad 
\hbox{for $u,v\in \cH_\pi$.}
$$
We note that 
$$
\norm{\pi(F)}\le \norm{F}_{L^1}\quad\hbox{and}\quad \pi(F)^\ast=\pi(F^\ast),
$$
where $F^\ast(g):=\overline{F(g^{-1})}$. 
The map $F\mapsto \pi(F)$ is a $\ast$-representation of the convolution algebra $L^1(G)$ into the algebra $\cB(\cH_\pi)$ of bounded operators on $\cH_\pi$.  

Let $\Gamma$ be a lattice subgroup in $G$. 
We use the Haar measure $m_{\Gamma \backslash G}$ on $\Gamma \backslash G$ such that
\begin{equation}\label{eq:measures}
\int_{\Gamma \backslash G}\left(\sum_{\gamma\in\Gamma} \phi(\gamma g)\right)\, dm_{\Gamma \backslash G}(\Gamma g)= \int_G \phi\, dm_G,\quad \phi\in C_c(G)\,.
\end{equation}
We normalize $m_G$, so that 
$m_{\Gamma \backslash G}$ is a probability measure.
Let us consider the unitary representation of $G$ in $L^2(\Gamma\setminus G)$, denoted $\lambda_{\Gamma\backslash G}$.
When the lattice $\Gamma$ is uniform, $\lambda_{\Gamma\backslash G}$ decomposes to a countable direct sum of irreducible sub-representations, each with finite multiplicity
$m(\pi,\Gamma)$ (cf. \cite{GGPS}): 
\begin{equation}\label{disc-decomp}
\lambda_{\Gamma\backslash G}=\bigoplus_{\pi\in \widehat{G}}m(\pi,\Gamma)\cdot \pi\,.
\end{equation}
Given $F\in L^1(G)$, the operator $\lambda_{\Gamma\setminus G}(F)$ is given by
\begin{equation}\label{eq:av}
\lambda_{\Gamma\backslash G}(F)(\phi)(\Gamma x)=\int_G F(g)\phi(\Gamma xg)\,dm_G(g),\quad \hbox{for $\phi\in L^2(\Gamma\backslash G)$ and $\Gamma x\in \Gamma\backslash G$.}
\end{equation}
We recall that when $F\in C_c(G)$,  the operator
$\lambda_{\Gamma\backslash G}(F)$
is an integral operator. Namely, for the continuous kernel
$$
L_F(\Gamma x, \Gamma y):=\sum_{\gamma\in \Gamma} F(x^{-1}\gamma y), \quad x,y\in G,
$$
we have
$$
\lambda_{\Gamma\backslash G}(F)\phi(\Gamma x)=\int_{ \Gamma \backslash G} L_F(\Gamma x, \Gamma y)\phi(\Gamma y)\,dm_{\Gamma\setminus G}(\Gamma y),
\quad \hbox{for $\phi\in L^2(\Gamma\backslash G)$ and $\Gamma x,\Gamma y\in \Gamma\backslash G$.}
$$
Therefore, for a function of the form $F=f\ast f^\ast$ with $f \in C_c(G)$, the operator 
$\lambda_{\Gamma\setminus G}(F)$ is a self-adjoint non-negative operator which is of trace class, and its trace is given by (cf. \cite{GGPS}):
$$\text{Tr}\big(\lambda_{\Gamma\setminus G}(F)\big)=\int_{\Gamma \setminus G}
L_F(\Gamma x, \Gamma x)\,dm_{\Gamma\setminus G}(\Gamma x)= \int_{ \Gamma \setminus G}\left(\sum_{\gamma \in \Gamma} F(x^{-1} \gamma  x)\right)\, dm_{\Gamma\setminus G}(\Gamma x).
$$
Combining this with (\ref{disc-decomp}), the {\it pre-trace formula} follows: 
\begin{equation}\label{pre-trace}
\text{Tr}\big(\lambda_{\Gamma\setminus G}(F)\big)=\sum_{\pi \in \widehat{G}} m(\pi, \Gamma) \text{Tr}(\pi(F))=\int_{\Gamma \setminus G}\left(\sum_{\gamma \in \Gamma} F(x^{-1} \gamma x)\right)\, dm_{\Gamma\setminus G}(\Gamma x)\,.
\end{equation}

Let $K$ be a compact subgroup of $G$ such that 
$(G, K)$ is a Gelfand pair.
An irreducible representation $\pi\in \widehat{G}$ is called
spherical if there exists a unit vector $v_\pi\in \cH_\pi$
which is invariant under $K$.
This vector is unique up to a complex phase.
We denote by $\widehat{G}^{sph}\subset \widehat{G}$
the subset  of spherical representations.
Given $\pi\in \widehat{G}^{sph}$ and a $K$-bi-invariant function $F\in L^1(G)$, 
$v_\pi$ is an eigenvector of the operator $\pi(F)$. 
The eigenvalue is given in terms of the spherical function 
$\varphi_\pi(g):=\inn{v_\pi, \pi(g)v_\pi}$ by the formula 
\begin{equation}\label{sph-eignevalue}
\pi(F)(v_\pi)=\varphi_\pi(F)\cdot v_\pi\,,  \text{   where  } \varphi_\pi(F):=\int_G F(g)\overline{\varphi_\pi(g)}\,dm_G(g). 
\end{equation}
The map $F\mapsto \varphi_\pi(F)$ is a $\ast$-homomorphism of $L^1(G,K)$, the convolution algebra of $K$-bi-invariant functions on $L^1(G)$.
In particular,
\begin{equation}\label{eq:starr}
\int_G (F^\ast\ast F) (g)\cdot  \overline{\varphi_\pi(g)}\, dm_G(g)=\abs{\int_G F(g)\overline { \varphi_\pi(g)}\, dm_G(g)}^2\,.    
\end{equation}
For $F\in L^1(G,K)$, the trace is given by
\begin{equation}\label{eq:trs}
\text{Tr}\big(\lambda_{\Gamma\setminus G}(F)\big)=\sum_{\pi \in \widehat{G}^{sph}} m(\pi, \Gamma) \text{Tr}( \pi(F))
=\sum_{\pi \in \widehat{G}^{sph}} m(\pi, \Gamma) \varphi_\pi(F)\,.
\end{equation}

Now we suppose that $G=G_1 \times G_2$ and $K=K_1\times K_2$,
so that $(G_1, K_2)$ and $(G_1, K_2)$
are Gelfand pairs.
Then every irreducible representation $\pi\in \widehat{G}$
is of the form $\pi=\pi_1\otimes \pi_2$
with $\pi_1\in \widehat{G}_1$ and $\pi_2\in \widehat{G}_2$.
Furthermore, $\pi$ is $K$-spherical if and only if $\pi_i$ is $K_i$-spherical for $i=1,2$.
Then the spherical function $\varphi_\pi$ has the form 
$$
\varphi_\pi(g_1, g_2)=\varphi_{\pi_1}(g_1)\varphi_{\pi_2}(g_2), \quad (g_1,g_2)\in G_1\times G_2.
$$
For $F=f_1\otimes f_2$ with $f_1\in L^1(G_1,K_1)$ and
$f_2\in L^1(G_2,K_2)$,
\begin{equation}\label{prod-pos-def}  
\varphi_\pi(F^\ast \ast F)
=\abs{\varphi_\pi(F)}^2= 
\abs{\varphi_{\pi_1}\left( f_1\right)}^2\cdot \abs{\varphi_{\pi_2}\left( f_2\right)}^2.
\end{equation}

\subsection{General spectral condition for the optimal exponent}\label{sec:general}

Let $G_i$, $i=1,2$, be lcsc non-compact groups, 
and $K_i\subset G_i$
are compact subgroups
such that $(G_i, K_i)$ are Gelfand pairs. 
Consider
the corresponding homogeneous spaces $X_i:=G_i/K_i$, which we equip with proper 
$G_i$-invariant 
metrics $D_i$. 
Given an irreducible lattice $\Gamma$ in $G=G_1\times G_2$, we consider the action of $\Gamma$ on the space
$X_1$. 
For  given points $x_1,y_1\in G_1$, 
we will study the existence of solutions $\gamma=(\gamma_1,\gamma_2)\in\Gamma$
for the inequalities
\begin{equation}\label{eq:inn}
D_1(\gamma_1 x_1K_1,y_1K_1)\le  r \quad\hbox{and}\quad \quad D_2(\gamma_2 K_2,K_2) \le R
\end{equation}
with the parameters $r\to 0$ and $R\to \infty$ at  appropriate rates.

Our analysis is based on investigation of
the operators $\lambda_{\Gamma\backslash G}(F_{r,R})$,
where $F_{r,R}:=f_{1,r}\otimes f_{2,R}$
with the functions $f_{1,r}$ and $f_{2,R}$ being   
smooth approximations to the bi-$K_i$-invariant sets
\begin{equation}\label{eq:ball}\widetilde{B}^{G_1}_r:=\set{g_1\in G_1\,;\, D_1(K_1, g_1K_1) \le  r}\quad,\quad \widetilde{B}^{G_2}_R:=\set{g_2\in G_2\,;\, D_2(K_2, g_2K_2) \le R}.
\end{equation}
Our first main result provides a solution to \eqref{eq:inn}
once a suitable bound on the Hilbert--Schmidt norm of the operators
$\lambda_{\Gamma\backslash G}(F_{r,R})$ is established. Keeping the notation and assumptions just introduced, we can now state : 

\begin{theorem}\label{th:general}
Let $\Gamma$ be a uniform irreducible lattice in $G=G_1\times G_2$, $R\in (1,\infty)$,
and  let ${\bf r}:(1,\infty)\to (0,1)$ be a non-increasing function.
We assume that
\begin{enumerate}
    \item[(SG)] The action of $G_2$ on $L^2_0(\Gamma\backslash G)$ has the following spectral gap property: there exist 
    $\mathfrak{w}>0$ and non-negative  
 $K_2$-bi-invariant functions $\psi_{2,R}\in C_c(G_2)$ satisfying 
    $$
\supp(\psi_{2,R})\subset\widetilde{B}^{G_2}_{R},\quad  
\int_{G_2} \psi_{2,R}\, dm_{G_2}=1\,,
\quad
\big\|\lambda_{\Gamma\backslash G}(\psi_{2,R})|_{L_0^2(\Gamma\backslash G)}\big\| \le const\cdot e^{-\mathfrak{w} R}    $$
where the implied constant is independent of $R$. 
    \item[(T)]  For every $\delta > 0$,   there exist 
    $(K_1\times K_2)$-invariant functions $F_{{\bf r}(R), R}\in C_c(G)$ such that 
$$
F_{{\bf r}(R), R}:=f_{1,{\bf r}(R)}\otimes f_{2,R}\quad\hbox{ with $f_{1,{\bf r}(R)}\in  C_c(G_1)$ and $f_{2,R}\in C_c(G_2)$}
$$
satisfying
$$
\supp(f_{1,{\bf r}(R)})\subset \widetilde B^{G_1}_{{\bf r}(R)},\quad \supp(f_{2,R})\subset \widetilde B^{G_2}_{R},\quad \int_G F_{{\bf r}(R), R}\, dm_G=1,
$$
and 
$$
\hbox{\rm Tr}\Big(\lambda_{\Gamma\backslash G}(F^*_{{\bf r}(R), R}*F_{{\bf r}(R),R})\Big) \ll_\delta e^{c\delta R}
$$
for  all sufficiently large $R$, with $c$ independent of $\delta$ and $R$.
\end{enumerate}
Then for every $\delta>0$ and almost every $(x_1,y_1)\in G_1\times G_1$,
the inequalities 
$$
D_1(\gamma_1 x_1K_1,y_1 K_1)\le  2{\bf r}(R) ,\quad D_2(\gamma_2 K_2,K_2) \le (1+\delta)R
$$
have solutions $\gamma=(\gamma_1,\gamma_2)\in\Gamma$ for all $R>R_0(x_1,y_1,\delta)$.
\end{theorem}

We remark that below we will construct functions $f_{1,{\bf r}(R)}\in  C_c(G_1)$ and $f_{2,R}\in C_c(G_2)$ such that
$$
\supp(f_{1,{\bf r}(R)})\subset \widetilde B^{G_1}_{{\bf r}(R)}\quad\hbox{and}\quad  \int_{G_1} f_{1,{\bf r}(R)} \, dm_{G_1}\asymp m_{G_1}\big(\widetilde B^{G_1}_{{\bf r}(R)}\big)^{1/2},
$$
and 
$$
\supp(f_{2,R})\subset \widetilde B^{G_2}_{R}\quad\hbox{and}\quad  \int_{G_2} f_{2,R} \, dm_{G_2}\asymp_\delta m_{G_2}\big(\widetilde B^{G_2}_{R}\big)^{1/2\pm\delta}\;\;\hbox{for all $\delta>0$.}
$$
Then the condition $\int_G F_{{\bf r}(R),R}\, dm_G=1$ 
will be reduced to the {\it matching volumes condition}:
$$
m_{G_1}\big(\widetilde B^{G_1}_{{\bf r}(R)}\big)\cdot m_{G_2}\big(\widetilde B^{G_2}_{R}\big)\asymp 1.
$$
Since 
$$
m_{G_1}\big(\widetilde B^{G_1}_{r}\big)\asymp r^{\dim(X_1)}\;\;
\hbox{for $r\in (0,1)$}
\quad
\hbox{and}
\quad
m_{G_2}\big(\widetilde B^{G_2}_{R}\big)\asymp e^{(a(G_2)\pm \delta)R}\;\;
\hbox{for $R\in (1,\infty)$,}
$$
and the matching volume condition will correspond to the estimate \eqref{eq:oopt}.

\subsection{Optimal approximation exponents for quaternion lattices}

Let us now turn to state our second main result which establishes the optimal approximation exponent  for irreducible uniform arithmetic congruence lattices $\Gamma$  in 
$$
G:=\prod_{j=1}^\ell G^{(j)}\quad\hbox{with $G^{(j)}$ isomorphic either to $\hbox{SL}_2(\RR)$ or to $\hbox{SL}_2(\CC)$.}
$$
We refer to Section \ref{sec:quaternions} for a discussion of these lattices. Respectively, let
$$
X:=\prod_{j=1}^\ell X^{(j)}\quad\hbox{with $X^{(j)}$ either 
$\mathbb{H}^2$ or $\mathbb{H}^3$,}
$$
so that $G$ acts on $X$ by isometries.
We denote by $D^{(j)}$
the unique constant-curvature $-1$ metrics on the hyperbolic spaces $X^{(j)}$. We denote by $D$ the metric on $X_{[1,\ell]}$ given by 
\begin{equation}
\label{metricD}
{D}(x,y):=\max \Big(D^{(j)}(x_1^{(j)},y_1^{(j)}):\, j=1,\ldots, l\Big)\quad  \hbox{ for  
$x_1,y_1\in X_{[1,\ell]},$}
\end{equation}
We let $\cD$ denote the unique left $G$-invariant and $K$-bi-invariant metric on $G$ satisfying that 
$\max_{k\in K}\cD(g,g^\prime k)=D(gK, g^\prime K)$, and then $\abs{D(gK, g^\prime K)- \cD(g,g^\prime)}\le C$ uniformly.

Fix $k=1,\ldots, \ell-1$ and consider the factors
$$
X_{[1,k]}:=\prod_{j=1}^k X^{(j)}\quad\hbox{and}\quad X_{[k+1,\ell]}:=\prod_{j=k+1}^\ell X^{(j)},
$$
so that $X=X_{[1,k]}\times X_{[k+1,\ell]}$. 

 $\Gamma$ naturally acts on $X_{[1,k]}$, and 
we establish the following optimal density for this action:

\begin{theorem}\label{th:quaternion}
Let
$$
d_{[1,k]}:=\sum_{j=1}^k \dim_\RR (X^{(j)})\quad\hbox{and}\quad a_{[k+1,\ell]}:=\sum_{j=k+1}^\ell \big(\dim_\RR (X^{(j)})-1\big), 
$$
and  fix reference points $x_0^{(j)}\in X^{(j)}$, $j=k+1,\ldots,\ell$.
Then for every $\zeta> \frac{d_{[1,k]}}{a_{[k+1,\ell]}}$ and almost every $x=(x^{(1)},\ldots, x^{(k)}),y=(y^{(1)},\ldots, y^{(k)})\in {X}_{[1,k]}$, the inequalities
\begin{align*}
&D^{(j)}(\gamma^{(j)}\cdot x^{(j)},y^{(j)})\le \epsilon ,\;\; j=1,\ldots,k,\\
&D^{(j)}(\gamma^{(j)}\cdot x_0^{(j)},x_0^{(j)})\le \zeta \, \log(1/\epsilon), \quad j=k+1,\ldots,\ell,
\end{align*}
have solutions $\gamma=(\gamma^{(1)},\ldots,\gamma^{(\ell)})\in\Gamma$ for all $\epsilon\in (0,  \epsilon_0(x,y,\zeta))$. 

Furthermore, the approximation exponent for the $\Gamma$-action on $X_{[1,k]}$ with respect to the metic $\cD$ is optimal, namely, 
for almost every $x,y\in {X}_{[1,k]}$,
$$
\kappa_\Gamma(x,y)=\frac{d_{[1,k]}}{a_{[k+1,\ell]}}.
$$
\end{theorem}

\begin{remark}
We note that all irreducible arithmetic lattices are commensurable to those arising from quaternion algebras over number fields (see \cite[10.3.6,10.3.7]{McLR} or \cite[3.3]{Bor81}). Furthermore, the ones that are uniform correspond to the non-split quaternion algebras \cite[3.3]{Bor81}, i.e., ones not isomorphic to the matrix algebra. We discuss these lattices in Section \ref{sec:quaternions} below.
\end{remark}


\section{Proof of Theorem \ref{th:general}}
\label{sec:proof_general}

\subsection{Geometric estimates}
\label{sec:geom}
For $r,R>0$ and a compact subset $C_1\subset  G_1$, we introduce the sets
\begin{equation}\label{eq:oomega}
\Omega_{r,R}(C_1):=\big\{ (x_1,y_1)\in C_1\times C_1:\,\, \not\exists \gamma\in \Gamma:\, D_1(\gamma_1 x_1 K_1, y_1 K_1) \le r,\,  D_2(\gamma_2 K_2,K_2)\le R\big\}.
\end{equation}
Our goal in the present section is to give a quantitative estimate on the measures of these sets as the parameters $r\to 0$ and $R\to \infty$ in a prescribed fashion.

Let $\psi_{1,r}$ be a non-negative continuous $K_1$-bi-invariant function on $G_1$ such that 
$$
\supp(\psi_{1,r})\subset\widetilde{B}^{G_1}_{r}\quad\hbox{and}\quad \int_{G_1} \psi_{1,r}\, dm_{G_1}=1.
$$

Similarly, for a fixed $R_0 > 0$, we introduce a non-negative continuous $K_2$-bi-invariant function $\psi_{2,R_0}$ on $G_2$ such that 
$$
\supp(\psi_{2,R_0})\subset\widetilde{B}^{G_2}_{R_0}\quad\hbox{and}\quad \int_{G_2} \psi_{2,R_0} \, dm_{G_2}=1.
$$
Since the parameter $R_0$ will be fixed throughout, we omit it from our notations.
Let
$$
\psi_r(g_1,g_2)=\psi_{1,r}(g_1)\psi_{2,R_0}(g_2)\quad \hbox{for $(g_1,g_2)\in G_1\times G_2.$}
$$
Further, let $x=(x_1, x_2)\in G$ and $y=(y_1,y_2)\in G$. We set
\begin{equation}\label{eq:phii}
\phi_{x,r}(y) :=\sum_{\gamma\in\Gamma} \psi_r(x^{-1}\gamma^{-1} y)=\sum_{(\gamma_1,\gamma_2)\in\Gamma} \psi_r(x_1^{-1}\gamma^{-1}_1y_1, x_2^{-1}\gamma^{-1}_2y_2).
\end{equation}
This defines a continuous compactly supported function on 
the homogeneous space $Z:=\Gamma\backslash G$
which is $(K_1\times K_2)$-invariant under the right action of $G$. Below we will also view it as a left-$\Gamma$-invariant function on $G$. 
We note that 
 \begin{equation}\label{normalize} 
 \int_G \psi_r \, dm_G=1 \text{     and so         }\int_{Z} \phi_{x,r}\, dm_{Z}=1
\end{equation} 
due to our normalization \eqref{eq:measures}.
 
The following lemma allows us to reduce the question about existence of solutions
for inequalities with $\gamma\in \Gamma$ to analysis of suitable averaging 
operators of the form \eqref{eq:av}.

\begin{lemma}\label{lem:1}
Fix parameters $r,R_0>0$, $r',R>0$, and $L>0$.
Let $x=(x_1,x_2)\in G$ and $y=(y_1,y_2)\in G$ such that $x_2,y_2\in \widetilde B^{G_2}_{L}$.
Let $f=f_{1}\otimes f_2\in C_c(G)$ such that 
$$
\supp(f_1)\subset \widetilde{B}^{G_1}_{r^\prime}\quad\hbox{and}\quad \supp(f_2)\subset \widetilde{B}^{G_2}_{R}.
$$
Then provided that
$$
\lambda_{Z}(f)\phi_{x,r}(y)\ne 0,
$$
it follows that 
there exists $\gamma=(\gamma_1,\gamma_2)\in\Gamma$ such that 
$$
D_1(\gamma_1 x_1K_1, y_1K_1)\le r+r' \quad\hbox{and}\quad D_2(\gamma_2 K_2,K_2)\le  2L+R_0+R.
$$
\end{lemma}

\begin{proof}
We observe that our assumption implies that 
$$
\int_{G_1\times G_2} \left(\sum_{(\gamma_1,\gamma_2)\in\Gamma} \psi_{1,r}(x_1^{-1} \gamma^{-1}_1 y_1 g_1)\psi_{2,R_0}(x_2^{-1}\gamma^{-1}_2 y_2 g_2) \right)f_1(g_1)f_2(g_2)\, dm_{G_1}(g_1) dm_{G_2}(g_2)\ne 0.
$$
Therefore, there exist $g_1\in G_1$ and $g_2\in G_2$ 
satisfying 
$$
D_1(g_1K_1,K_1)\le r'\quad\hbox{and}\quad D_2(g_2K_2,K_2)\le R
$$
such that 
$$
D_1(x_1^{-1} \gamma^{-1}_1 y_1 g_1 K_1,K_1)\le r\quad\hbox{and} \quad D_2(x_2^{-1}\gamma^{-1}_2 y_2 g_2 K_2,K_2)\le R_0
$$
for some $(\gamma_1,\gamma_2)\in \Gamma$.
Using invariance of the metric, it follows from the first inequality that
$$
D_1(\gamma_1 x_1K_1,y_1K_1)\le D_1(\gamma_1 x_1K_1,y_1g_1K_1)+D_1(y_1g_1K_1,y_1K_1)\le r+r'.
$$
Also it follows from the second inequality that
\begin{align*}
D_2(\gamma_2 K_2,K_2)&\le D_2(\gamma_2 K_2,y_2g_2K_2) +D_2(y_2g_2K_2,K_2)\\
&\le \big(D_2(K_2,x_2K_2)+D_2(x_2 K_2,\gamma_2^{-1}y_2g_2K_2)\big)+
\big(D_2(g_2K_2, K_2)+D_2(K_2,y_2^{-1} K_2)\big) \\
&\le 2L+R_0+R.
\end{align*}
This completes the proof of the lemma.
\end{proof}

Taking Lemma \ref{lem:1} into account, we can now estimate
the measure of the sets \eqref{eq:oomega} in terms of
the averaging operators:

\begin{lemma}\label{lem:2}
Let $r,r',R>0$, $f\in C_c(G)$ and $\phi_{x,r}\in C_c(G)$ be as in Lemma
\ref{lem:1}, and $C_1\subset G_1$ and $C_2\subset G_2$ are 
compact subsets.
Then there exists $c_0>0$, depending only on $C_2$ and $R_0$, such that 
$$
m_{G_1\times G_1}\big(\Omega_{r+r^\prime,R+c_0}(C_1)\big)\ll_{C_1,C_2} \int_{C_1\times C_2} \big\|\lambda_{Z}(f)\phi_{x,r}-1\big\|^2_{L^2(Z)}\, dm_{G}(x).
$$
\end{lemma}

\begin{proof}
For $r,R>0$ and $x_1\in C_1$, we write
$$
\Omega^{x_1}_{r,R}(C_1):=\big\{ y_1\in C_1:\,\, \not\exists \gamma\in \Gamma:\, D_1(\gamma_1 x_1 K_1, y_1 K_1)\le r,\,  D_2(\gamma_2 K_2,K_2)\le R\big\}.
$$
Then 
$$
m_{G_1\times G_1}(\Omega_{r,R}(C_1)))
=\int_{x_1\in C_1} m_{G_1}(\Omega^{x_1}_{r,R}(C_1))\, dm_{G_1}(x_1).
$$
Therefore, it is sufficient to show that for all $x=(x_1,x_2)\in C_1\times C_2$, 
\begin{equation}\label{eq:sss}
m_{G_1}\big(\Omega^{x_1}_{r+r^\prime,R+c_0}(C_1)\big)\ll_{C_1,C_2}  \big\|\lambda_{Z}(f)\phi_{x,r}-1\big\|^2_{L^2(Z)},
\end{equation}
which we proceed to verify.  

Suppose that 
$$
\lambda_{Z}(f)\phi_{x,r}(y)\ne 0
$$
for some  $x=(x_1,x_2)\in C_1\times C_2$ and $y=(y_1,y_2)\in C_1\times C_2$.
Then it follows from Lemma \ref{lem:1} that 
there exists $\gamma=(\gamma_1, \gamma_2)\in\Gamma$ such that 
$$
D_1(\gamma_1 x_1K_1, y_1K_1)\le r+r' \quad\hbox{and}\quad D_2(\gamma_2 K_2,K_2)\le 2L+R_0+R,
$$
where the constant $L$ is determined by the subset $C_2$.
Let us set $c_0=2L+R_0$. 
We conclude that 
for every $y_1\in \Omega^{x_1}_{r+r^\prime,R+c_0}(C_1)$
and $y_2\in C_2$,
$$
\lambda_{Z}(f)\phi_{x,r}(y_1,y_2)= 0.
$$
Therefore, it follows that
$$
m_{G_1}\big(\Omega^{x_1}_{r+r^\prime,R+c_0}(C_1)\big) m_{G_2}(C_2)\le 
\big\|\lambda_{Z}(f)\phi_{x,r}-1\big\|^2_{L^2(C_1\times C_2)}\,.
$$
Finally, we observe that for every $\phi\in C_c(Z)$,
$$
\big\| \phi\circ pr\big\|_{L^2(C_1\times C_2)}\ll_{C_1,C_2} \big\|\phi \big\|_{L^2(Z)},
$$
where $pr:G\to Z=\Gamma\backslash G:g\mapsto \Gamma g$ denotes the natural factor map.
Therefore, (\ref{eq:sss}) follows directly from the last estimate.
\end{proof}

\subsection{Norm bounds}
\label{sec:spec}

In this section we also use the functions $\phi_{x,r}\in C_c(Z)$ introduced 
in \eqref{eq:phii}. Recall that they are normalized so that 
$\int_Z \phi_{x,r}\, dm_Z=1$.

\begin{lemma}\label{lem:3}
Let $r,r'>0$ and $R\ge 1$. Fix positive constants $\theta$ and $\alpha$, and consider the following two claims:

\begin{enumerate}
    \item[(I)] For every $\delta>0$, there exists a 
    $(K_1\times K_2)$-bi-invariant function $F=f_1\otimes f_2\in C_c(G)$ such that 
    $$
    \supp(f_1)\subset \widetilde{B}^{G_1}_{r'},\quad \supp(f_2)\subset \widetilde{B}^{G_2}_{R},\quad \int_G F\, dm_G=1,
    $$
    and for
    compact subsets $C_1\subset G_1$, $C_2\subset G_2$, and all $R\ge 1$,
$$
\int_{C_1\times C_2}\big\|\lambda_{Z}(F)\phi_{x,r}\big\|^2_{L^2(Z)}\, dm_{G}(x) \ll_{C_1,C_2,\delta} e^{\delta R}.
$$

\item[(II)] There exists 
a 
$(K_1\times K_2)$-bi-invariant function
$ F^\prime=h_1\otimes h_2\in C_c(G)$ such that 
$$
\supp(h_1)\subset \widetilde{B}^{G_1}_{r'},\quad
\supp(h_2)\subset \widetilde{B}^{G_2}_{(1+\theta)R},\quad  \int_G F^\prime\, dm_G=1,
$$
and for compact subsets $C_1\subset G_1$, $C_2\subset G_2$, and all $R\ge 1$,
$$
\int_{C_1\times C_2}\big\|\lambda_{Z}(F^\prime)\phi_{x,r}-1\big\|^2_{L^2(Z)}\, dm_{G}(x) \ll_{C_1,C_2,\theta} e^{- \alpha R}. 
$$

\end{enumerate} 
Then for every $\theta > 0$ there exists $\alpha > 0$ such that Claim I implies Claim II. 
\end{lemma}

\begin{proof}
First, we observe that (I) implies that
\begin{align*}
&\;\;\;\;\;\int_{C_1\times C_2}\big\|\lambda_{Z}(F)\phi_{x,r}-1\big\|^2_{L^2(Z)}\, dm_{G}(x) \\
&\le
\int_{C_1\times C_2}\left(\big\|\lambda_{Z}(F)\phi_{x,r}\big\|^2_{L^2(Z)}+2\big\|\lambda_{Z}(F)\phi_{x,r}\big\|_{L^2(Z)}+1\right)\, dm_{G}(x) \\
&\le 
\int_{C_1\times C_2}\big\|\lambda_{Z}(F)\phi_{x,r}\big\|^2_{L^2(Z)}\, dm_{G}(x) \\
&\;\;\;\;\;+ 2
\left(\int_{C_1\times C_2}\big\|\lambda_{Z}(F)\phi_{x,r}\big\|^2_{L^2(Z)}\, dm_{G}(x)\right)^{1/2} m_{G}(C_1\times C_2)^{1/2}+
m_{G}(C_1\times C_2)\\
&\ll_{C_1,C_2,\delta} e^{\delta R}.
\end{align*}
Let us introduce an additional parameter $R^\prime>0$
and choose a non-negative $K_2$-bi-invariant 
function $\psi_{2, R^\prime}$ on $G_2$ such that 

$$
\supp(\psi_{2,R^\prime})\subset\widetilde{B}^{G_2}_{R^\prime},\quad  
\int_{G_2} \psi_{2,R^\prime} \, dm_{G_2}=1,\quad
\big\|\lambda_{Z}(\psi_{2,R^\prime})|_{L_0^2(Z)}\big\| \le const\cdot e^{-\mathfrak{w} R^\prime}
$$
for some fixed $\mathfrak{w}>0$. The existence of $\psi_{2,R^\prime}$ follows from our assumption (SG), namely, the spectral gap 
property for the action of $G_2$ on $L^2_0(Z)$.
Let $F^\prime=f_1\otimes (\psi_{2,R^\prime}\ast f_2)$. Then  
$$
\supp(\psi_{2,R^\prime}\ast f_2)\subset \widetilde{B}^{G_2}_{R+R^\prime}\quad\hbox{and}\quad \int_{G} F^\prime \, dm_G=1\,,
$$
and furthermore, since $\psi_{2,R^\prime}\ast (f_1\otimes f_2)=f_1\otimes (\psi_{2,R^\prime}\ast f_2)$, we have 
$$\lambda_{Z}(F^\prime)= \lambda_Z(\psi_{2,R^\prime})\lambda_Z(F).$$
Because of our normalization conditions on $F$, $F^\prime$, $\psi_{2,R^\prime}$ and $\phi_{x,r}$,
$$
\lambda_{Z}(F^\prime)\phi_{x,r}-1\in L^2_0(Z)\,\text{  and  } \lambda_Z(F^\prime)1=1\,.
$$
We conclude that
\begin{align*}
\int_{C_1\times C_2}\big\|\lambda_{Z}(F^\prime)\phi_{x,r}-1\big\|^2_{L^2(Z)}\, dm_{G}(x) &=
\int_{C_1\times C_2}\big\|\lambda_{Z}(\psi_{2,R^\prime})(\lambda_{Z}(F)\phi_{x,r}-1)\big\|^2_{L^2(Z)}\, dm_{G}(x) \\
&\le const\cdot  e^{-2\mathfrak{w} R^\prime} 
\int_{C_1\times C_2}\big\|\lambda_{Z}(F)\phi_{x,r}-1\big\|^2_{L^2(Z)}\, dm_{G}(x) \\
&\ll_{C_1,C_2,\delta} e^{-2\mathfrak{w} R^\prime}\cdot e^{\delta R}. 
\end{align*}
Now set $R^\prime=\theta R$ and $\delta=\mathfrak{w}\theta$.
This gives the required bound in (II) with $\alpha=\mathfrak{w}\theta> 0$.
\end{proof}

\medskip

From now on, we assume that $\Gamma$ is a uniform lattice, and as usual, that $(G_1,K_1)$ and $(G_2, K_2)$ are Gelfand pairs. Then the space
$L^2(Z)$ decomposes as a direct sum of  irreducible representations $\pi=\pi_1\otimes \pi_2$ of $G_1\times G_2$
with finite multiplicities $m(\pi,\Gamma)$. Let $\widehat G^{\rm sph}$ denote the spherical unitary dual of $G$, consisting of irreducible spherical representations, namely those containing a unit vector invariant under $K_1\times K_2$. 
For $\pi\in \widehat G^{\rm sph}$, we denote by $v_\pi$ the unique (up to a complex phase) $(K_1\times K_2)$-invariant unit vector.
Since the function $\phi_{x,r}$, defined in \eqref{eq:phii}, is $(K_1\times K_2)$-invariant under the right action of $G$ on $Z$, 
it can be written as the following $L^2$-convergent series. 
Let $\lambda_Z=\bigoplus m(\pi, \Gamma)\pi$ be the decomposition of $\lambda_Z$
into isotypical components $\cH_\pi$ of irreducible representations. We retain only those that are spherical, denoting this set by $\Sigma_{\Gamma\setminus G}^{sph}$. Then the $\pi$-isotypical component $\cH_\pi$ contains $m(\pi, \Gamma)$ mutually orthogonal unit vectors, each of which is an eigenvector of every  convolution operator $\pi(F)$, with $F\in L^1(G,K)$. We denote these eignevectors in $\cH_\pi$ by $v_\pi^{(j)}$, $j=1,\dots, m(\pi,\Gamma)$. We now have the following orthogonal expansion 
\begin{equation}\label{eq:sph}
\phi_{x,r}=\sum_{\pi\in \Sigma_{\Gamma\setminus G}^{sph}}\sum_{j=1}^{m(\pi,\Gamma)} \left<\phi_{x,r}, v_{\pi^{(j)}} \right> v_{\pi^{(j)}},
\end{equation}

This will allow us to compute explicitly the action of the averaging operators on $\phi_{x,r}$, as follows.

\begin{lemma}\label{l:spectral}
Let $F=f_1\otimes f_2$, where $f_i$ are $K_i$-bi-invariant functions with 
compact support on $G_i$, $i=1,2$. Then 
for any compact $C_1\subset G_1$ and $C_2\subset G_2$,
$$
\int_{C_1\times C_2} \big\|\lambda_{Z}(F)\phi_{x,r}\big\|^2_{L^2(Z)}\, dm_{G}(x)\ll_{C_1,C_2} \sum_{\pi\in \Sigma_{\Gamma\setminus G}^{sph}} m(\pi,\Gamma)
\abs{\left<\pi(F)v_\pi,v_\pi\right>}^2
$$
uniformly for $r>0$.
\end{lemma}

\begin{proof}
It follows from the orthogonality of the expansion \eqref{eq:sph} that 
\begin{equation}\label{sum-squares}
\big\|\lambda_{Z}(F)\phi_{x,r}\big\|^2_{L^2(Z)}=\left<\lambda_{Z}(F)\phi_{x,r}, \lambda_{Z}(F)\phi_{x,r} \right>
=
\sum_{\pi\in \Sigma_{\Gamma\setminus G}^{sph}} m(\pi,\Gamma)\abs{\left<\phi_{x,r}, v_\pi \right>}^2 
\left<\pi(F)v_\pi, \pi(F)v_\pi \right>.
\end{equation}
Since the space of $(K_1\times K_2)$-invariant vectors for a spherical irreducible representation $\pi$ is one-dimensional, 
\begin{equation}\label{component} 
\pi(F)v_\pi=\left<\pi(F)v_\pi,v_\pi\right>v_\pi,
\end{equation}
so that
\begin{equation}\label{sum-squares2}
\big\|\lambda_{Z}(F)\phi_{x,r}\big\|^2_{L^2(Z)}=
\sum_{\pi\in \Sigma_{\Gamma\setminus G}^{sph}} m(\pi,\Gamma)\abs{\left<\phi_{x,r}, v_\pi \right>}^2 
\abs{\left<\pi(F)v_\pi,v_\pi\right>}^2.
\end{equation}
 
 Using the definition of $\phi_{x,r}$ (see \eqref{eq:phii}), we write (replacing $\gamma^{-1}$ by $\gamma$) 
$$
\left<\phi_{x,r},v_\pi\right>=\int_{\Gamma\backslash G}  \left(
\sum_{\gamma\in\Gamma} \psi_r(x^{-1}\gamma y)\overline{v_\pi(y)}\right)dm_{\Gamma\backslash G}(\Gamma y),
$$
and now viewing $v_\pi$ as a left $\Gamma$-invariant function on $G$ we have 
$$\left<\phi_{x,r},v_\pi\right>=
\int_G \psi_r(x^{-1}g)\overline{v_\pi(g)}\, dm_G(g)
= 
\int_G \psi_r(g)\overline{v_\pi(xg)}\, dm_G(g).
$$
We recall that for $\Gamma y\in\Gamma\backslash G$,
$$
\pi(\psi_r)v_\pi(\Gamma y)=\int_G \psi_r(g)v_\pi(yg)dm_G(g),
$$
so that (since $\psi_r$ is real)
$$
\left<\phi_{e,r},v_\pi\right>=\overline{\pi(\psi_r)v_\pi(e)}.
$$
It follows from the uniqueness of $(K_1\times K_2)$-invariant vectors $v_\pi$ that 
$$
\pi(\psi_r)v_\pi=\left<\pi(\psi_r)v_\pi,v_\pi\right>v_\pi.
$$
Hence, we conclude that 
$$
\left<\phi_{e,r},v_\pi\right>=\overline{\left<\pi(\psi_r)v_\pi,v_\pi\right>}\cdot\overline{v_\pi(e)}.
$$
This provides a formula for $\left<\phi_{e,r},v_\pi\right>$.

For general $x\in G$, consider the map $J_x : \Gamma\backslash G\to \left(x^{-1}\Gamma x\right)\backslash G$ 
given by $J_x(\Gamma g)=\left(x^{-1}\Gamma x\right)x^{-1}g$, 
which is a well defined $G$-equivariant isomorphism of $G$-spaces.
Now define 
$$
\Pi_x: L^2(\Gamma\backslash G)\to L^2(\left(x^{-1}\Gamma x\right)\backslash G): \phi\mapsto \phi \circ J_x^{-1},
$$
which is an isomorphism of unitary representations, between $\pi$ and the unitary representation of $G$ on $L^2(\left(x^{-1}\Gamma x\right)\backslash G)$, which we denote by $\pi_x$. 
Then $ \left(\Pi_x F\right)(x^{-1}\Gamma x\cdot  g)=F(\Gamma xg)$, for $F\in L^2(\Gamma \backslash G)$. 
We set
$$
\widetilde \phi_{x,r}(y)=\sum_{\gamma\in\Gamma} \psi_r(x^{-1}\gamma x y)\quad\hbox{and}\quad
\widetilde v_\pi(y)=v_\pi(xy).
$$
Here we view $\widetilde v_\pi(y)$ as a function on $G$, and note that it is $x^{-1}\Gamma x$-invariant under left translations, since 
$$\widetilde v_\pi(x^{-1}\gamma x y)=v_\pi(x(x^{-1}\gamma xy))=v_\pi(\gamma xy)=v_\pi(xy)=\widetilde v_\pi(y).$$ 
In addition, $\widetilde v_\pi(y)$ is invariant under right ranslations by $K_1\times K_2$ since $v_\pi$ has this property. 
It follows that  $\widetilde v_\pi$ is (the lift to $G$ of) the unique  $(K_1\times K_2)$-invariant unit vector
for the representation $\pi_x$.
Therefore, it follows from the previous argument  that 
\begin{align*}
\left<\phi_{x,r},v_\pi\right>=\left<\Pi_x(\phi_{x,r}),\Pi_x(v_\pi)\right>
=\left<\widetilde\phi_{x,r},\widetilde v_\pi\right>
=\overline{\left<\pi_x(\psi_r)\widetilde v_\pi,\widetilde v_\pi\right>}\cdot\overline{\widetilde v_\pi(e)}
=\overline{\left<\pi(\psi_r) v_\pi, v_\pi\right>}\cdot\overline{ v_\pi(x)}.
\end{align*}
Using the last formula, we can rewrite (\ref{sum-squares2}) as 
$$
\big\|\lambda_{Z}(F)\phi_{x,r}\big\|^2_{L^2(Z)}
=
\sum_{\pi\in \Sigma_{\Gamma\setminus G}^{sph}} m(\pi,\Gamma)\abs{\left<\pi(\psi_r) v_\pi, v_\pi \right>}^2 |v_\pi(x)|^2
\abs{\left<\pi(F)v_\pi,v_\pi\right>}^2.
$$
Now since the functions $v_\pi$ is $L^2$-normalized, 
$$
\int_{C_1\times C_2} |v_\pi(x)|^2\, dm_{G}(x)\ll_{C_1,C_2}
\int_{Z} |v_\pi|^2\, dm_{Z}=1,
$$
and since additionally $\|\psi_r\|_{L^1(G)}=\int_G \psi_r\, dm_G=1$,
$$
\big|\left<\pi(\psi_r) v_\pi, v_\pi \right>\big|\le \|\pi(\psi_r)\|\le 1. 
$$
Therefore, we conclude that 
$$
\int_{C_1\times C_2} \big\|\lambda_{\Gamma\setminus G}(F)\phi_{x,r}\big\|^2_{L^2(Z)}\, dm_{G}(x)\ll_{C_1,C_2} \sum_{\pi\in \Sigma_{\Gamma\setminus G}^{sph}} m(\pi,\Gamma) 
\abs{\left<\pi(F)v_\pi,v_\pi\right>}^2,
$$
and this completes the proof of the lemma.
\end{proof}

\subsection{Completion of the proof of Theorem \ref{th:general}}
\label{sec:proof0}

Let $R>1$ and $\delta > 0$.
According to our assumption (T), there exist 
$(K_1\times K_2)$-invariant functions $F_{{\bf r}(R),R}\in C_c(G)$ such that 
$F_{{\bf r}(R),R}=f_{1,{\bf r}(R)}\otimes f_{2,R}$ with $f_{1,{\bf r}(R)}\in  C_c(G_1)$ and $f_{2,R}\in C_c(G_2)$ satisfying
$$
\supp(f_{1,{\bf r}(R)})\subset \widetilde B^{G_1}_{{\bf r}(R)},\quad \supp(f_{2,R})\subset \widetilde B^{G_2}_{R},\quad \int_G F_{{\bf r}(R),R}\, dm_G=1,
$$
and 
\begin{equation}\label{eq:ttt}
\hbox{Tr}\Big(\lambda_Z(F_{{\bf r}(R),R}^\ast\ast F_{{\bf r}(R),R})\Big)=\sum_{\pi\in \widehat G^{\rm sph}} m(\pi,\Gamma)\left<\pi(F_{{\bf r}(R),R}^\ast\ast F_{{\bf r}(R), R})v_\pi,v_\pi\right>\ll_\delta e^{\delta R}\,. 
\end{equation}
for all sufficiently large $R$.
We also note that by (\ref{component}),
$$ 
\inn{\pi(F_{{\bf r}(R),R}^\ast\ast F_{{\bf r}(R), R})v_\pi, v_\pi}=\inn{\pi(F_{{\bf r}(R), R})v_\pi, \pi(F_{{\bf r}(R), R})v_\pi}=\abs{\left<\pi(F_{{\bf r}(R), R})v_\pi,v_\pi\right>}^2.
$$

Let us fix compact subsets $C_1\subset G_1$ and $C_2\subset G_1$. Then 
it follows from Lemma \ref{l:spectral}  and \eqref{eq:ttt} that
$$
\int_{C_1\times C_2} \big\|\lambda_{Z}(F_{{\bf r}(R), R})\phi_{x,r}\big\|^2_{L^2(Z)}\, dm_{G}(x)\ll_{C_1,C_2,\delta}
e^{\delta R}\,. 
$$
Since $\delta>0$ is arbitrary, this verifies Property (I) of 
Lemma \ref{lem:3}.
Next, we take an arbitrary $\theta>0$.
By Lemma \ref{lem:3} with $r'=r={\bf r}(R)$, there exists $\alpha=\alpha(\theta)>0$ and functions $F^\prime_{{\bf r}(R), R}=h_{1,{\bf r}(R)}\otimes  h_{2,R}\in C_c(G)$ such that 
$$
\supp(h_{{\bf r}(R)})\subset \widetilde{B}_{1,{\bf r}(R)}\quad\hbox{and}\quad \supp( h_{2,R})\subset \widetilde{B}^{G_2}_{(1+\theta)R}, $$
and 
$$
\int_{C_1\times C_2}\big\|\lambda_{Z}(F^\prime_{{\bf r}(R),R})\phi_{x,r}-1\big\|^2_{L^2(Z)}\, dm_{G}(x) \ll_{C_1,C_2,\theta} e^{-\alpha R}. 
$$
Hence, we deduce from Lemma \ref{lem:2} that
$$
m_{G_1\times G_1}\big(\Omega_{2{\bf r}(R),(1+\theta)R+c_0}(C_1)\big)\ll_{C_1} e^{-\alpha R}
$$
for sufficiently large $R$. 
Once we have this bound, it remains to apply a Borel-Cantelli argument.

We consider the sets 
$$
\Omega_N:=\Omega_{2{\bf r}(N),(1+\theta)N+c_0}(C_1)
$$
for $N\in \mathbb{N}$. It follows from Borel--Cantelli Lemma that for almost every $(x_1,y_1)\in C_1\times C_1$,
we have $(x_1,y_1)\notin \Omega_N$ for all $N>N_0(x_1,y_1,\theta)$, namely, 
the inequalities 
$$
D_1(\gamma_1 x_1K_1,y_1K_1)\le  2{\bf r}(N) ,\quad D_2(\gamma_2 K_2,K_2) \le (1+\theta)N+c_0
$$
have solutions $\gamma=(\gamma_1,\gamma_2)\in\Gamma$.
For general $R>1$, we pick $N\in\mathbb{N}$ such that $N-1\le R<N$.
Then, when $R$ is sufficiently large,
the inequalities 
$$
D_1(\gamma_1 x_1K_1,y_1K_1)\le  2{\bf r}(N)\le 2{\bf r}(R)  ,\quad D_2(\gamma_2 K_2,K_2) \le (1+\theta)N+c_0\le 
(1+\theta)(R+1)+c_0
$$
have solutions $\gamma=(\gamma_1,\gamma_2)\in\Gamma$.
Since $\theta$ is an arbitrary positive parameter, we conclude that 
for every $\delta>0$ and almost every $(x_1,y_1)\in C_1\times C_1$,
the inequalities 
$$
D_1(\gamma_1 x_1K_1,y_1K_1)\le  2{\bf r}(R) ,\quad D_2(\gamma_2 K_2,K_2) \le (1+\delta)R
$$
have solutions $\gamma=(\gamma_1,\gamma_2)\in\Gamma$ when $R>R_0(x_1,y_1,\delta)$.
Finally, since the compact set $C_1$ is arbitrary, this completes the proof of the Theorem  \ref{th:general}. \qed

\section{Quaternion lattices}\label{sec:quaternions} 

Let $m,n\in \NN_{\ge 0}$ be natural numbers with $m+n> 0$. In this section we give a construction of arithmetic lattices in $G=\SL_2(\RR)^m\times \SL_2(\CC)^n$ following the classical paper by Borel \cite{Bor81}. Before that, let us recall what does it mean for a general lattice in a real semisimple group $G$ to be \emph{arithmetic}. 
\begin{definition}\label{def-arithmetic}
    Let $G$ be a conected semisimplie Lie group with finite center. A lattice $\Gamma\subset G$ is called arithmetic if there exists a semisimple algebraic group $\mathbb G$ defined over $\mathbb Q$, together with a $\QQ$-embedding $\mathbb G\hookrightarrow \GL_n$ and a projection $\text{proj}\colon \mathbb G(\mathbb R)\to G$, such that 
    \begin{enumerate}
        \item $\text{proj}(\mathbb G(\mathbb Z))$ is commensurable with $\Gamma$,
        \item the kernel of the projection  is compact.
    \end{enumerate}
\end{definition}

In addition, we say that a lattice is \emph{irreducible}, if the projection of $\Gamma$ onto every proper factor of $G$ is dense. In the case $G=\SL_2(\RR)^m\times \SL_2(\CC)^n$ and irreducible $\Gamma$, all groups $\mathbb G$ that can appear in the above definition arise from quaternion algebras defined over number fields (see \cite[\S 8]{McLR} and \cite{Bor81}). Let us turn to review some elements of the construction. 
\begin{definition}
Let $\FF$ be a field of characteristic zero. A quaternion algebra $\cA$ over $\FF$ is the algebra determined by the relations: 

$$\cA=\left(\frac{\mathfrak{u},\mathfrak{v}}{\FF}\right):=\FF+{\bf i}\FF+{\bf j}\FF+{\bf ij}\FF, \quad {\bf i}^2=\mathfrak{u}, {\bf j}^2=\mathfrak{v}, {\bf ij}=-{\bf ji},$$ where $\mathfrak{u},\mathfrak{v}\in \FF$ are non-zero.
\end{definition}

The simplest example is the $2\times 2$ matrix algebras. Indeed, 
$$ {\rm M}_2(\FF)\simeq \left(\frac{1,1}{\FF}\right),\text{ with } {\bf i}=\begin{pmatrix}1 & 0\\ 0 & -1\end{pmatrix} \text{ and } {\bf j}=\begin{pmatrix}0 & 1\\ 1 & 0\end{pmatrix}.$$
In general ${\bf i},{\bf j}$ in $\left(\frac{\mathfrak{u},\mathfrak{v}}{\FF}\right)$ can be chosen as  $2\times 2$-matrices over a quadratic extension of $\FF$ :
$${\bf i}=\begin{pmatrix} \sqrt{\mathfrak{u} } & 0\\ 0 & -\sqrt{\mathfrak{u}} 
\end{pmatrix},\quad {\bf j}=\begin{pmatrix}
0 & 1\\ \mathfrak{v} & 0   
\end{pmatrix}.$$

A quaternion algebra $\cA:=\left(\frac{\mathfrak{u},\mathfrak{v}}{\FF}\right)$ is isomorphic to ${\rm M}_2(\FF)$ if and only if the equation $x_1^2-\mathfrak{u}x_2^2-\mathfrak{v}x_3^2+\mathfrak{u}\mathfrak{v}x_4^2=0$ admits a non-trivial solution over $\FF$. In that case, we say that $\cA$ splits. If $\cA$ does not split, then it is a division algebra over $\FF$. In this case, we say that $\cA$ is not ramified over $\FF$. Same notations apply to extensions of $\FF$. 

Any quaternion algebra over $\FF$ becomes isomorphic to ${\rm M}_2(\overline \FF)$, over the algebraic closure $\overline \FF$ of $\FF$. The determinant and trace on ${\rm M}_2(\overline \FF)$ restrict respectively to the norm ${\bf n}\colon \cA\to \FF$ and trace ${\bf tr}\colon \cA\to \FF$ on $\cA$.
$$ {\bf tr}(x+{\bf i}y + {\bf j}z +{\bf ij}t):= 2x,\quad  {\bf n}(x+{\bf i}y + {\bf j}z +{\bf ij}t):=x^2-\mathfrak{u}y^2-\mathfrak{v}z^2+\mathfrak{u}\mathfrak{v} t^2.$$
Both functions ${\bf tr}, {\bf n}$ are rational over $\FF$ and are additive and multiplicative respectively.

Let us now choose $\FF$ to  be a number field $\KK$ of degree $\mathfrak{d}$ and let $\mathfrak o_\KK$ be the ring of integers of $\KK$.  The \emph{ramification set} of $\cA$ is the set ${\rm ram}(\cA)$ of all equivalence classes of valuations $|\cdot|_v$ on $\KK$ such that $\cA_v:=\cA\otimes_\KK \KK_v$ does not split. The ramification set determines the isomorphism class of a quaternion algebra \cite{McLR}. For example ${\sf H}:=\left(\frac{-1,-1}{\QQ}\right)$ is a non-split quaternion algebra over $\QQ$. The ramification set of $\sf H$ consists of the unique real place and a single finite place corresponding to the prime $2$. Over the unique real place ${\sf H}\otimes_\QQ \RR=\mathcal H=\left(\frac{-1,-1}{\RR}\right)$ is the familiar algebra of Hamilton quaternions. 

In order to use $\cA$ to construct lattices in $\SL_2(\RR)^m\times \SL_2(\CC)^n$ we have to impose certain conditions on the ramification of $\cA$ in the Archimedean places of $\KK$. We say that $\cA$ is \emph{admissible} if $\KK$ has exactly $m$ real Archimedean valuations where $\cA$ splits and $n$ complex Archimedean valuations where $\cA$ splits. The remaining $\mathfrak{d}-m-2n$ real valuations of $\KK$ are in ${\rm ram}(\cA)$. This means that 
$$\cA\otimes_\QQ \mathbb R\simeq \prod_{v}\cA_v={\rm M}_2(\RR)^m \times {\rm M}_2(\CC)^n\times \mathcal H^{\mathfrak{d}-m-2n},$$
where $\mathcal H$ is the algebra of Hamilton's quaternions. 

The norm function allows us to define the group $$\mathbb G=\SL_1(\cA):=\{ x\in \cA\mid {\bf n}(x)=1\}.$$ It is a simple algebraic group defined over $\KK$ which becomes isomorphic to $\SL_2(\overline \KK)$ over the algebraic closure of $\KK$. Strictly speaking, Definition \ref{def-arithmetic} requires a group defined over $\QQ$ which in this case is obtained by a restriction of scalars from $\KK$. Our construction will bypass the restriction of scalars but still produce arithmetic lattices in the sense of Definition \ref{def-arithmetic}. The splitting conditions satisfied by an admissible quaternion algebra $\cA$ guarantee that 
$$\mathbb G(\KK\otimes _\QQ \RR)\simeq \SL_2(\RR)^m\times \SL_2(\CC)^n\times \SU(2)^{\mathfrak{d}-m-2n}.$$

Let $\text{proj}\colon \mathbb G(\KK\otimes _\QQ \RR)\to \SL_2(\RR)^m\times \SL_2(\CC)^n$ be the projection onto the first $m+n$ factors. We note that the kernel of $\text{proj}$ is $\SU(2)^{\mathfrak{d}-m-2n}$, so it is compact. 

The last ingredient is to define an arithmetic subgroup of $\mathbb G(\KK\otimes _\QQ \RR)$ which will play the role of integral points from the Definition \ref{def-arithmetic}. Instead of embedding $\mathbb G$ into $\GL_n$ we will do this using orders in quaternion algebras. An \emph{order} in $\cA$ is a unital $\mathfrak o_\KK$-subalgebra $\mathfrak O$ of $\cA$ such that $\mathfrak O\otimes_{\mathfrak o_\KK} \KK=\cA.$ A classical example is provided by the order of Hurwitz quaternions
$$\mathfrak H:=\{ x+ {\bf i}y +{\bf j}z+{\bf ij}t \mid x,y,z,t\in \mathbb Z\text{ or } x,y,z,t\in \mathbb Z+{1}/{2}\} \subset \left(\frac{-1,-1}{\QQ}\right).$$

Given an order $\mathfrak O$ of $\cA$ we consider the group of units 
$$\mathfrak O^1:=\{ x\in \mathfrak O\mid {\bf n}(x)=1\}.$$
Since $\mathfrak O^1\subset \mathbb G(\KK)$, we can consider it as subgroup of $\mathbb G(\KK\otimes_\QQ \RR)$. By the Borel-Harish-Chandra theorem \cite{BHC62}, the group $\mathfrak O^1$ is a lattice in $\mathbb G(\KK\otimes _\QQ \RR)$. The kernel of the projection $\text{proj}\colon \mathbb G(\KK\otimes_\QQ \RR)\to \SL_2(\RR)^m\times \SL_2(\CC)^n$ is compact, so $\mathfrak O^1$ projects to a lattice in $\SL_2(\RR)^m\times \SL_2(\CC)^n$, denoted $\Gamma_{\mathfrak O}$: 
$$\Gamma_{\mathfrak O}:=\text{proj}(\mathfrak O^1)\subset \SL_2(\RR)^m\times \SL_2(\CC)^n.$$
The lattice $\Gamma_{\mathfrak O}$ is cocompact if and only if $\cA$ is a division algebra over $\KK$ \cite[Theorem 8.12]{McLR}. 

\begin{definition}An irreducible lattice $\Lambda\subset \SL_2(\RR)^m\times \SL_2(\CC)^n$ is called \emph{arithmetic} if it commensurable to some $\Gamma_{\mathfrak O}$, possibly after conjugation.
\end{definition}

\begin{example}\label{ex-Qlattice}
Let $\KK=\QQ(\sqrt{17})$. Then $\cA=\left(\frac{3,5}{\KK}\right)$ is a non-split quaternion algebra with ramification at $3,5$ (which both remain prime in $\QQ(\sqrt{17}))$. The algebra $\cA$ is admissible for $G=\SL_2(\RR)^2.$  The elements ${\bf i},{\bf j}$ can be chosen as 
$${\bf i}=\begin{pmatrix}
\sqrt{3} & 0 \\ 0 & -\sqrt{3}
\end{pmatrix}, \quad {\bf j}=\begin{pmatrix}
0 & 1 \\ 5 & 0
\end{pmatrix}.$$ Let $\omega:=\frac{1+\sqrt{17}}{2}$. The ring of integers of $\KK$ is $\mathfrak o_\KK=\ZZ[\omega]$. We choose an order $$\mathfrak O=\ZZ[\omega]+{\bf i}\ZZ[\omega]+{\bf j}\ZZ[\omega]+{\bf ij}\ZZ[\omega].$$ Then $\mathfrak O$ consists of 
\begin{align*}&\begin{pmatrix}
u_1 + v_1 \omega +u_2\sqrt{3}+v_2\sqrt{3}\omega & u_3+v_3\omega+u_4\sqrt{3}+v_4\sqrt{3}\omega\\
5u_3+5v_3\omega -5u_4\sqrt{3}-5v_4\sqrt{3}\omega & u_1 + v_1\omega -u_2\sqrt{3}-v_2\sqrt{3}\omega
\end{pmatrix}\end{align*}
with $u_1,u_2,u_3,u_4,v_1,v_2,v_3,v_4\in \mathbb Z$.
To obtain $\Gamma_{\mathfrak O}$, we embed $\mathfrak O^1$ diagonally in $\SL_2(\RR)\times \SL_2(\RR)$ via maps sending $\omega$ to $\frac{1+\sqrt{17}}{2}$ and $\frac{1-\sqrt{17}}{2}$ respectively.
\end{example}

Let $\mathfrak I$ be an ideal of $\mathfrak o_\KK$. Define 
$$\mathfrak O^1(\mathfrak I):=\{x\in \mathfrak O^1\mid x-1\in \mathfrak I\mathfrak O\}\text{ and } \Gamma_{\mathfrak O}(\mathfrak I):=\text{proj}(\mathfrak O^1(\mathfrak I)).$$ Taking the order $\mathfrak O$ as in Example \ref{ex-Qlattice} and an infinite family of ideals $\mathfrak I$ of $\mathbb Z[\omega]$ we get an infinite family of lattices to which Theorem \ref{th:quaternion} applies. 
\begin{definition} An irreducible lattice $\Lambda\subset \SL_2(\RR)^m\times \SL_2(\CC)^n$ is an \emph{arithmetic congruence lattice} if it contains $\Gamma_{\mathfrak O}(\mathfrak I),$ up to conjugation, for some ideal $\mathfrak I$. 
\end{definition}


\section{Proof of Theorem \ref{th:quaternion}}
\label{sec:proof_quat}

\subsection{Preliminaries}
In this section, we work with the groups $\hbox{SL}_2(\RR)$ and  $\hbox{SL}_2(\CC)$. To have unified notation, we introduce a parameter 
\begin{equation}\label{eq:rho}
\rho:=\left\{ 
\begin{tabular}{ll}
$\frac{1}{2}$, & \hbox{for $\hbox{SL}_2(\RR)$,}\\
$1$, & \hbox{for $\hbox{SL}_2(\CC)$.}
\end{tabular}
\right.
\end{equation}
(We note that this differs from the notation in \cite[\S 3]{FHMM},
but we prefer to follow this convention because the above $\rho$ 
corresponds as customary to the half-sum of positive roots.)
Respectively, $L_\rho$ denotes either 
$\hbox{SL}_2(\RR)$ or  $\hbox{SL}_2(\CC)$, and 
$K_\rho$ denotes either
$\hbox{SO}_2(\RR)$ or $\hbox{SU}_2(\CC)$.
Then $K_\rho$ is a maximal compact subgroup of $L_\rho$
and 
$$
A:=\set{a_t:=e^{tH_1}: t\in \RR},\,\,\text{with} \,\, H_1:=\left(\begin{smallmatrix}1/2 & 0\\
0 &-1/2\end{smallmatrix}\right),
$$ 
is a Cartan subgroup of $L_\rho$ compatible with $K_\rho$.
Let $\ma$ be the Lie algebra of $A$. 
We use the identification 
\begin{equation}\label{eq:ident}
\Lie{a}\simeq \RR:\, tH_1\mapsto t.
\end{equation}
We note that $H_1$ is the element of $\ma$ satisfying $\alpha(H_1)=1$ for the unique positive root $\alpha$. 

The groups $L_{\rho}$ act by isometries on the hyperbolic spaces $\HH^{2\rho+1}$
and $\hbox{Stab}_{L_{\rho}}(i)=K_{\rho}$, so that 
$$
\HH^{2\rho+1}\simeq L_{\rho}/K_{\rho}.
$$
We set 
$$
d_{\rho}:=2\rho+1=\dim_\RR (\HH^{2\rho+1}).
$$
We have the Cartan decomposition 
$$
L_{\rho}=K_{\rho}A K_{\rho},
$$
and a Haar measure on $L_{\rho}$ is given by: for $f\in C_c(G)$,
\begin{equation}\label{eq:cartan}
\int_{L_{\rho}}  f\, dm_{L_{\rho}} =c_\rho \int_{K_{\rho}}\int_\RR\int_{K_{\rho}} f(k a_t k')( \sinh t)^{2\rho}\,dm_{K_{\rho}}(k)\, dt  \,dm_{K_{\rho}}(k')\,,
\end{equation}
for some $c_\rho>0$. 

Let $D_\rho$ denote the standard hyperbolic metric on $\HH^{2\rho+1}$ and 
$$
\widetilde{B}^{L_{\rho}}_r:=\set{g\in L_{\rho} \,:\, D_\rho(g\cdot i,i) \le r}.
$$
We recall that $\{ka_t\cdot i=k\cdot e^ti:\, t\ge 0\}$, with $k\in K_{\rho}$, are
the unit speed geodesic rays starting at $i$, so that 
$$
\widetilde{B}^{L_{\rho}}_r=K_{\rho}\set{a_t\,;\, 0\le t \le r}K_{\rho}.
$$
In particular, it follows from \eqref{eq:cartan} that
\begin{equation}\label{volume-SL2} 
m_{L_{\rho}}(\widetilde{B}^{L_{\rho}}_r)\asymp \begin{cases}  e^{2\rho r}, & \text{ if } r\geq 1,\\
 r^{2\rho+1}, & \text{ if } 0\leq r\le 1.
\end{cases} 
\end{equation}

The unitary spherical dual of $L_\rho$ is parameterized by
$$
\widehat L_\rho^{sph}\simeq (0,1/2]\cup i[0,\infty).
$$
For $z$ in this set, we denote by
$\pi_z$ the corresponding irreducible spherical unitary representation of $L_\rho$. The point $z=1/2$ gives the trivial representation, the points
$z\in (0,1/2]$ give the non-tempered complimentary series,
and the points $z\in i[0,\infty)$ give the tempered
principal series.
We denote by $\varphi_z=\varphi_z^{(\rho)}$ the spherical function of the representation $\pi_z$.

For an irreducible non-trivial unitary representation $\pi$ of $L_{\rho}$, we denote by 
$p^+(\pi)$ the {\it integrability exponent} of $\pi$
$$
p^+(\pi):=\inf\set{2\le p \,;\, \hbox{$\pi$ has a non-zero   matrix coefficient in $L^p(L_{\rho})$} }\,.
$$
In the case when the representation of $\pi_z$ is spherical as above, the integrability exponent is given by 
$$p^+(\pi_z)=\inf\set{2\le p \,;\, \varphi_z \in L^p(L_{\rho})}=\frac{2}{1-2\re (z)}$$
(see e.g. \cite[\S 1.3.3 (1.11)]{FHMM}).

\subsection{Smooth approximations for balls} \label{sec:spherical functions}

Our identification \eqref{eq:ident} gives
$$
\ma^*_\CC:=\ma^*+i\ma^*\simeq \CC.
$$
For every $z\in \CC$, we denote by $\varphi_z^{(\rho)}$
the elementary spherical function on $L_{\rho}$ corresponding to the parameter $z$. 
It can be realized as the matrix coefficient of the spherical principal series of (not necessarily unitary) representation of $L_{\rho}$ induced from the character
$a_t\mapsto e^{2\rho z t}$, namely 
$\varphi_z^{(\rho)}(g)=\inn{v_z, \pi_z(g)v_z}$, where $v_z$ is the unique (up to complex phase) $K_{\rho}$-invariant unit vector in the induced representation. In this notation, the imaginary axis $z=i\tau\in i\RR$  corresponds to the tempered spectrum, while the real interval $z=s\in [-\frac12, \frac12] $ corresponds to the complementary series, see \cite[\S 3.3]{FHMM}. 
  In particular  $\varphi^{(\rho)}_{\frac12}=1$, and $\varphi^{(\rho)}_0$ is the Harish Chandra $\Xi$-function on $L_{\rho}$.  Note that $\varphi^{(\rho)}_z$ is well-defined even if $z$ does not correspond to any actual unitary representation.

Let $f$ be a compactly supported, bounded, $K_\rho$-bi-invariant function on $L_{\rho}$. The spherical transform of $f$ is a holomorphic function on $\CC$,  defined by 
\begin{equation}\label{spher-tran} \widehat{f}(z):= \int_{L_{\rho}} f(g) \varphi^{(\rho)}_z(g^{-1})\, dm_{L_{\rho}}(g).
\end{equation}

We recall the Paley-Wiener theorem \cite{Gangolli} for the spherical transform in the  version given in  \cite[\S 3.3]{FHMM}:

 \begin{theorem}\label{Paley-Wiener}
 \begin{enumerate}
 \item[(a)] Suppose that $f$ is a smooth $K_{\rho}$-bi-invariant function on $L_{\rho}$ 
 with 
 \begin{equation}\label{eq:support}
 \supp(f)\subset \widetilde{B}^{L_{\rho}}_R.
 \end{equation}
 Then the spherical transform $\widehat{f}^{}(z)$ is an entire holomorphic function of $z\in \CC$, symmetric under $z\mapsto -z$, and 
 for every $N>0$, there exists a constant $C_N>0$, depending on $f$, such that 
$$ \abs{\widehat{f}^{}(z)}\leq C_N\, e^{2\rho R\abs{\sigma}}(1+\abs{\tau})^{-N}\quad \hbox{for all $z=\sigma+i\tau \in \CC$.}
$$ 

\item[(b)]
Conversely, if $\varphi$ is an entire holomorphic function on $\CC$, invariant under $z \mapsto -z$ 
and for every $ N > 0$ there exist constants $C_N>0$ such that $$ \abs{\varphi(z)}\leq C_N\, e^{2\rho R\abs{\sigma}}(1+\abs{\tau})^{-N}
\quad \hbox{for all $z=\sigma+i\tau \in \CC$,}
$$ then $\varphi=\widehat{f}^{}$ is the spherical transform of a unique smooth $K_{\rho}$-bi-invariant function $f$ on $L_{\rho} $ 
satisfying \eqref{eq:support}.
\end{enumerate}
\end{theorem}
The constants $C_N$ in Theorem \ref{Paley-Wiener} depend on 
the function $f$, in general. We will apply these estimates below for a restricted family of functions $f$ defined geometrically, and for a restricted subset of $\CC$ (corresponding to the unitary spherical spectrum of $L_{\rho}$), where the constants will, in fact, be uniform. 

  We now turn to define two families of functions, a small scale family $f^{(\rho)}_{1,r}$, and a large scale family $f^{(\rho)}_{2,R}$, with precisely prescribed supports and spherical transform decay properties.   

\begin{lemma}\label{l:small}
Let $\delta^\prime\in (0,1)$. Then there exists a family of smooth $K_{\rho}$-bi-invariant functions $f^{(\rho)}_{1,r}$, $r\in (0,1)$, 
on $L_{\rho}$ such that 
\begin{equation}\label{integral-1}
\supp(f^{(\rho)}_{1,r})\subset \widetilde{B}^{L_{\rho}}_r,\quad
 \int_{L_{\rho}}  f^{(\rho)}_{1,r}\, dm_{L_{\rho}}\asymp_{\delta^\prime} r^{(2\rho+1)/2} \,,
\end{equation}
and for every $N>0$, there exists a constant $C_N>0$, independent of $r$, such that 
\begin{equation}\label{ball-estimator1}
\abs{\widehat{f}^{(\rho)}_{1,r}(z)}\le C_N\,  
r^{(2\rho+1)/2} (1+r\abs{\tau})^{-N}
\end{equation} 
for all $z=\sigma+i\tau\in \CC$ corresponding to the spherical unitary spectrum of $L_\rho$.
\end{lemma}

\begin{proof}
Let $\delta^\prime\in (0,1)$. We take a smooth even  function $\eta$ on $\RR$ such that 
$$
0\le \eta \le 1, \quad \supp(\eta)\subset [-1,1],\quad \eta=1\;\;\hbox{on $[-1+\delta^\prime,1-\delta^\prime]$}
$$
and set
$$
{f}^{(\rho)}_{1,r}(g):= r^{-(2\rho +1)/2}\cdot \eta \big(D_\rho(g\cdot i,i)/r\big),\quad r\in (0,1).
$$
Clearly, this is a smooth $K_\rho$-bi-invariant function
with support contained in $\widetilde{B}^{L_{\rho}}_r$.
Moreover, in terms of the Cartan decomposition,
$$
{f}^{(\rho)}_{1,r}(ka_t k^\prime)=r^{-(2\rho +1)/2}\cdot \eta (t/r),\quad \hbox{for $k,k^\prime\in K_\rho,\; t\in \RR^+.$}
$$
Therefore, it follows from \eqref{eq:cartan} that
$$
\int_{L_{\rho}} f^{(\rho)}_{1,r}\, dm_{L_{\rho}} =r^{-(2\rho +1)/2}\cdot \int_{\RR^+}\eta(t/r)(\sinh t)^{2\rho} \,dt $$
This implies the estimate on the integral of 
$f^{(\rho)}_{1,r}$.

Let us note that the family of functions $f^{(\rho)}_{1,r}$ depends on $\eta$ and hence on the parameter $\delta^\prime\in (0,1)$. 
Therefore the estimate of the integrals in (\ref{integral-1}) is up to two positive multiplicative constants depending on $\delta^\prime$, 
which converge to $1$ as $\delta^\prime \to 0$. The constants $C_N$ appearing in (\ref{ball-estimator1}) also depend on $\eta$ and $\delta^\prime$.
We will suppress this dependence in our notation, and return to discuss this point  later on in the proof of Theorem \ref{th:quaternion}.

Finally, we note that the estimate (\ref{ball-estimator1}) on the spherical transform of $f^{(\rho)}_{1,r}$ and its feature of uniformity, namely its independence of $r$, will be crucial to our argument. However, its full proof  involves elaborate analysis of the local spherical transform, and for this reason we will defer this task to Section \ref{sec:local SF}, where we execute the local analysis for all symmetric spaces. Here we only remark that the desired estimate (\ref{ball-estimator1})
follows from Theorem \ref{lem-1}.
Indeed, in the present case the dimension of the symmetric space of $L_\rho$ is equal to $2\rho+1$, and our analysis in Section \ref{sec:local SF} (see,
in particular, \eqref{omega}) applies to the
functions $\omega_r\in C_c(L_\rho)$ defined by
$$
\omega_r(k a_t k^\prime):= r^{-(2\rho+1)} \eta(t/r),\quad \hbox{for $k,k^\prime\in K$ and $t\in \RR^+$.}
$$
Since $f_{1,r}^{(\rho)}=r^{(2\rho+1)/2}\cdot \omega_r$,
this implies the required estimate.
\end{proof}

\begin{lemma}\label{l:large}
Let $\delta^\prime\in (0,1)$. 
There exists a family of smooth $K_{\rho}$-bi-invariant functions $f^{(\rho)}_{2,R}$, $R\ge 1$, 
on $L_{\rho}$ such that 
\begin{equation} \label{integral-2}
 \supp(f^{(\rho)}_{2,R})\subset \widetilde{B}^{L_{\rho}}_R,\quad
e^{(1-\delta^\prime)\rho R} \ll_{\delta^\prime} \int_{L_{\rho}}  f^{(\rho)}_{2,R}\, dm_{L_{\rho}}\ll e^{\rho R} \,,
\end{equation}
and for every $N>0$, there exists a constant $C_N>0$, independent of $R$, such that 
\begin{equation}\label{ball-estimator}
\abs{\widehat{f}_{2,R}^{(\rho)}(z)}\le C_N\,  e^{2\rho R\abs{\sigma}} 
\left(1+R\abs{\tau}\right)^{-N}\quad \hbox{for all $z=\sigma+i\tau \in \CC$}.
\end{equation} 
\end{lemma}

\begin{proof}

Let $\delta^\prime\in (0,1)$. We consider again the smooth even  function $\eta$ on $\RR$ such that 
$$
0\le \eta \le 1, \quad \supp(\eta)\subset [-1,1],\quad \eta=1\;\;\hbox{on $[-1+\delta^\prime,1-\delta^\prime]$}
$$
and consider its Fourier-Laplace transform :
$$
\widehat{\eta} (z):=\int_\RR \eta(t) e^{zt} dt\quad\hbox{for $z\in \CC$.}
$$
This is an even and entire holomorphic function.  
By the Paley-Wiener theorem for functions on $\RR$, for every $N > 0$
there exists a constant $C_N > 0$ such that 
$$\abs{\widehat{\eta} (z)}\le C_N\, e^{ \abs{\sigma}}
\left(1+\abs{\tau}\right)^{-N} \quad \hbox{for all $z=\sigma+i\tau \in \CC$}.
$$
This formulation of the Paley-Wiener theorem is due to H\"ormander. We refer to \cite[Ch.~I, \S 2, p. 15]{He84} for a relevant discussion. 

For all $R > 0$, the function $\eta_R(t):=\frac{1}{R}\eta(t/R)$ is supported on $[-R, R]$ and 
$$\widehat{\eta}_R(z)=\frac{1}{R}\int_\RR \eta(t/R)e^{ zt}dt =
\int_\RR \eta(t)e^{R zt}dt =
\,\widehat{\eta}(Rz),$$
so that for all $N>0$ and $z=\sigma+i\tau \in \CC$,
$$
\big|\widehat{\eta}_R (2\rho z)\big|\le C_N\,  e^{2\rho R\abs{\sigma}} 
\left(1+2\rho R\abs{\tau}\right)^{-N}.
$$
We emphasize that the constants $C_N$ here depend on $N$ and the function $\eta$, but not on the parameter $R$. Let us now assume that $R\ge 1$. 
Then the function $\widehat{\eta}_R (2\rho z)$ satisfies the conditions of Theorem \ref{Paley-Wiener}(b), so that 
there exists a smooth $K_{\rho}$-bi-invariant function $f^{(\rho)}_{2,R}$
such that its spherical transform is given by
$$
\widehat{f}^{(\rho)}_{2,R}(z)=\widehat{\eta}_R(2\rho z)\quad \hbox{for $z\in\CC$,}
$$
and its support $\supp (f^{(\rho)}_{2,R})$ is contained in the set $\widetilde{B}^{L_{\rho}}_R$.

Since $\varphi_{1/2}^{(\rho)}=1$, we have 
\begin{align}\label{f-integral}
\int_{L_{\rho}}  f^{(\rho)}_{2,R}\,dm_{L_{\rho}}
=\widehat{f}^{(\rho)}_{2,R}(1/2)=\widehat{\eta}_{R}(\rho)=\int_{-1}^1 \eta(t)e^{\rho Rt}\,dt. 
\end{align} 
This implies the required estimate on the integral (\ref{integral-2}) of $f^{(\rho)}_{2,R}$. 

Note that as in the previous Lemma, the family of functions $ f^{(\rho)}_{2,R}$ depends on $\eta$ and $\delta^\prime$, and so do the constants $C_N$. We will suppress this dependence in the notation, and return to discuss this point later on in the proof of Theorem \ref{th:quaternion}. 
\end{proof}

\subsection{Spectral gap and density bounds}\label{sec:densityestimate}

Let us now turn to consider the group
$$
G:=\prod_{j=1}^\ell G^{(j)} \quad\; \hbox{with  $G^{(j)}=\SL_2(\RR)$ or $G^{(j)}=\SL_2(\CC)$}.
$$
Since every irreducible unitary representation of $G$
is a tensor product $\pi^{(1)}\otimes \cdots\otimes\pi^{(\ell)}$ of
irreducible unitary representations $\pi^{(j)}$ of the factors,
the unitary spherical dual of $G$ is parameterized by
$$
\widehat G^{sph}\simeq \Big((0,1/2]\cup i[0,\infty)\Big)^{\ell}.
$$
For ${\bf z}=(z_1, \dots , z_{\ell})$ from this set,
we denote by $\pi_{\bf z}=\pi_{z_1}\otimes \dots \otimes \pi_{z_\ell}$ the corresponding unitary representation
of $G$. 
The integrability exponent of $\pi_{\bf z}$ is given by 
\begin{equation}\label{integrability}
p^+(\pi_{\bf z})=\frac{2}{1-2m({\bf z})}\quad \hbox{with}\quad m({\bf z}):=\max_{1\le j \le \ell} \re(z_j)\,.
\end{equation} 

Let $\Gamma$ be a uniform irreducible arithmetic congruence lattice in $G$. 
The unitary representation of $G$ in $L^2(\Gamma\setminus G)$, denoted $\lambda_{\Gamma\setminus G}$,  decomposes as a countable direct sum of irreducible sub-representations with finite multiplicities $m(\pi,\Gamma)$:
$$
\lambda_{\Gamma\setminus G}=\bigoplus_{\pi\in \widehat{G}}\; m(\pi,\Gamma)\cdot \pi\,.
$$
An important consideration in our analysis is the fact that 
the components of all non-trivial representations
appearing in this decomposition are uniformly isolated
from the trivial representation:

\begin{theorem}\label{strong SP}
There exists $\mathfrak{w}^\prime>0$ 
such that for every non-trivial irreducible spherical unitary representation $\pi_{\bf z}=\pi_{z_1}\otimes \cdots \otimes \pi_{z_\ell}$ appearing in the decomposition of $L^2(\Gamma\backslash G)$,
$$
\re(z_i)\le 1/2-\mathfrak{w}^\prime\quad\hbox{for all $i=1,\dots, \ell$.}
$$ 
\end{theorem}

This is a particular case of Property ($\tau$) \cite[Th. 3.1]{C03} --- the strong spectral gap for authomorphic representations. In our setting, this follows
from the Jacquet--Langlands correspondence
\cite{JL70} and the spectral gap of $\hbox{SL}_2$
\cite{BB11}.


Another crucial ingredient of our analysis is the density estimates
established in \cite{FHMM}:

\begin{theorem}[\cite{FHMM}, Th.~1 and Remark~1]\label{density hypothesis}
For 
a  subset $Q$ of $\set{1,\ldots ,\ell}$,
consider ${\boldsymbol \sigma}={\boldsymbol \sigma}(Q^c)=(\sigma_j)_{j\in Q^c}\in(0,1/2]^{\abs{Q^c}}$, and $T>0$, define
the following bounded subset of the spherical spectrum of $G$ : 
$$
\cB_Q({\boldsymbol \sigma}, T):=\set{\pi_{\bf z}\,;\, 
{\bf  z}\in {\prod}_{j\in Q^c} ( \sigma_j, 1/2]\times {\prod}_{j \in Q} i[0, T]}\,.
$$
Then for every $\delta > 0$, there exists $C_\delta > 0$ such that
\begin{equation}\label{eq-DensityTwo}
\sum_{\pi\in \cB_Q({\boldsymbol \sigma}, T)}  \mm(\pi,\Gamma)\le C_\delta\, {\prod}_{j \in  Q} (1+T)^{d_j(1-2m({\boldsymbol \sigma}))+\delta},
\end{equation}
where $m({\boldsymbol  \sigma}):=\max_{j\in Q^c}\sigma_j$, and
$d_j=2$ for $G^{(j)}=\hbox{\rm SL}_2(\RR)$ and
$d_j=3$ for $G^{(j)}=\hbox{\rm SL}_2(\CC)$.
\end{theorem}

\subsection{Main trace estimate}\label{sec:trace}

We now turn to formulate and prove our main spectral estimate. 
We keep notation from the previous sections. Additionally, $K^{(j)}$ denote the maximal compact subgroup of $G^{(j)}$
and $\rho_j$ is defined as \eqref{eq:rho}.

\begin{theorem}\label{trace-estimate}
Let $1\le k<\ell$,  $r\in (0,1)$, $R>1$, and $\delta\in (0,1)$.

For $1\le j\le k$, let $f_r^{(j)}$ be a smooth $K^{(j)}$-bi-invariant function on $G^{(j)}$ 
whose spherical transform satisfies:
for every $N>0$, there exists a constant $C_N>0$, independent of $r$, such that 
\begin{equation}\label{ball-estimator00}
\abs{\widehat{f}_{r}^{(j)}({\bf z})}
\le C_N\,  
r^{(2\rho_j+1)/2} (1+r\abs{\boldsymbol \tau})^{-N}
\quad \hbox{for all ${\bf z}={\boldsymbol \sigma}+i{\boldsymbol  \tau}\in\widehat{G}^{sph}$}.
\end{equation}
For $k+1\le j\le \ell$, let $f_R^{(j)}$ be a smooth $K^{(j)}$-bi-invariant function on $G^{(j)}$ 
whose spherical transform satisfies:
for every $N>0$, there exists a constant $C_N>0$, independent of $R$, such that 
\begin{equation}\label{ball-estimator0}
\abs{\widehat{f}_{R}^{(j)}(\bf z)}\le C_N\,  e^{2\rho_j R\abs{\boldsymbol \sigma}} 
\left(1+R\abs{\boldsymbol \tau}\right)^{-N}\quad \hbox{for all ${\bf z}={\boldsymbol \sigma}+i{\boldsymbol\tau}\in\widehat{G}^{sph}$}.
\end{equation}
Let
$$
F_{r,R}:=f_r^{(1)}\otimes\cdots \otimes f_r^{(k)} \otimes f_R^{(k+1)}\otimes\cdots \otimes f_R^{(\ell)}.
$$
Suppose that the parameters $R$ and $r$ satisfy 
\begin{equation}\label{eq:balance}
\prod_{j=1}^k m_{G^{(j)}}(\widetilde{B}^{(j)}_r)\ll_\delta e^{\delta R}\cdot \prod_{j=k+1}^\ell m_{G^{(j)}}(\widetilde{B}^{(j)}_R)^{-1}.
\end{equation}
Then
$$
\text{\rm Tr}\Big(\lambda_{\Gamma\setminus G}(F^*_{r,R}*F_{r,R})\Big)\ll_\delta e^{c\delta R}
$$
for a constant $c>0$, independent of $r,R$ and $\delta$.
\end{theorem}
Note that (\ref{eq:balance}) is a weaker form of the balanced volumes condition. 
\begin{proof}
{\bf Part I : The spectral estimator}. Let $d_j:=2\rho_j+1$, namely, $d_j$  denotes
the dimension of the symmetric space of $G^{(j)}$.
In view of \eqref{volume-SL2}, the assumption \eqref{eq:balance} 
is equivalent to
\begin{equation}\label{eq:rr}
 \prod_{j=1}^k r^{d_j} \ll_\delta e^{\delta R}\cdot \prod_{j=k+1}^\ell e^{-(d_j-1) R}.
\end{equation}
We set 
\begin{align*}
h_r^{(j)}:=(f_r^{(j)})^* * f_r^{(j)}\quad \hbox{for $1\le j\le k$,}\quad\hbox{and}\quad
h_R^{(j)}:=(f_R^{(j)})^* * f_R^{(j)}\quad \hbox{for $k+1\le j\le \ell$.}
\end{align*}
Then
$$
(F^*_{r,R}*F_{r,R})(g)=\prod_{j=1}^k h^{(j)}_r(g_j)\prod_{j=k+1}^\ell h^{(j)}_R(g_j),\quad g=(g_1,\ldots, g_\ell)\in G,
$$
The representation $\lambda_{\Gamma\backslash G}$ splits as a direct sum of irreducible representations $\pi$ of $G$.
Moreover, each $\pi$ is of the form $\pi=\otimes_{j=1}^\ell \pi^{(j)}$, where  $\pi^{(j)}$ are irreducible representations of $G^{(j)}$. Since the lattice is irreducible, if $\pi$ is not trivial,
then each of its factors $\pi^{(j)}$ is also non-trivial.
We also note that since $F_{r,R}$ is invariant under the maximal 
compact subgroup of $G$ then $\pi(F_{r,R})=0$ for any irreducible non-spherical representation $\pi$, so that only spherical representations contribute to the trace estimate. Therefore,
in view of \eqref{eq:trs},
\begin{align*}
\text{Tr}\Big(\lambda_{\Gamma\setminus G}(&F^*_{r,R}*F_{r,R})\Big)\\
=&\sum_{\pi \in \widehat{G}^{sph}} m(\pi, \Gamma)\, \text{Tr}\, \pi(F^*_{r,R}*F_{r,R})
=\sum_{\pi \in \widehat{G}^{sph}} m(\pi, \Gamma) \, \varphi_\pi(F^*_{r,R}*F_{r,R})\\
=&\sum_{\otimes_{j=1}^\ell \pi^{(j)} \in  \widehat{G}^{sph}}
m\big(\otimes_{j=1}^\ell \pi^{(j)},\Gamma\big)
 \prod_{j=1}^k \varphi_{\pi^{(j)}} \big(h^{{(j)}}_r\big) \prod_{j=k+1}^\ell
 \varphi_{\pi^{(j)}} \big(h^{{(j)}}_R\big).
\end{align*}
Let us assume that the above representation $\pi=\otimes_{j=1}^\ell \pi^{(j)}$ is not trivial. Then each of the factors $\pi^{(j)}$
is either the complementary series representation 
$\pi_s$ with $s\in (0,1/2)$ or the principal series representation
$\pi_{it}$ with $t\in [0,\infty)$. This set 
is parameterized by the union of $\prod_{j\in P} (0,\infty) \times\prod_{j\in P^c} [0,1/2)$ 
where $P$ runs of over all subsets 
of $\{1,\ldots,\ell\}$.
From now on we fix such a subset $P$ and the corresponding partitions
$$
\{1,\ldots,k\}=I_{pr}\sqcup I_{comp}\quad\hbox{and}\quad
\{k+1,\ldots,\ell\}=J_{pr}\sqcup J_{comp}
$$
such that $P=I_{pr}\sqcup J_{pr}$.
Then our estimate reduces to bounding each of the sums
$$
\sum_{{\bf t}(P)\,,\, {\bf s}(P^c)} m\big((\otimes_{j\in P} \pi_{it_j})\otimes  (\otimes_{j\in P^c} \pi_{s_j}) ,\Gamma\big)
\prod_{j\in I_{pr}} \widehat{h}^{{(j)}}_r(it_j) 
\prod_{j\in I_{comp}} \widehat{h}^{{(j)}}_r(s_j) 
\prod_{j\in J_{pr}} \widehat{h}^{{(j)}}_R(it_j) 
\prod_{j\in J_{comp}} \widehat{h}^{{(j)}}_R(s_j),
$$
where the sum runs over ${\bf t}(P)=(t_j)_{j\in P}$ and 
${\bf s}(P^c)=(s_j)_{j\in P^c}$ such that the corresponding multiplicities
are not trivial. 

Since the spherical transform is $*$-homomorphism (cf. \eqref{eq:starr}),
$$
\widehat{h}^{{(j)}}_R=\abs{\widehat{f}^{{(j)}}_R}^2
\quad\hbox{and}\quad \widehat{h}^{{(j)}}_r=\abs{\widehat{f}^{{(j)}}_r}^2.
$$
Thefore, using (\ref{ball-estimator0}) we have the following estimates:
for $R> 1$, $t\in (0,\infty)$, and $s\in [0,1/2)$,\begin{align*}
\widehat{h}^{{(j)}}_R(it)\ll_{N} (1+Rt)^{-N}\;\hbox{ for all $N\in\NN$}\quad\hbox{and}\quad
\widehat{h}^{{(j)}}_R(s)\ll e^{4\rho_j Rs},   
\end{align*}
and using (\ref{ball-estimator00}) for $r\in (0,1)$, $t\in (0,\infty)$, and $s\in [0,1/2)$,
\begin{align*}
\widehat{h}^{{(j)}}_r(i t)\ll_{N} r^{2\rho_j+1}(1+rt)^{-N}\;\hbox{ for all $N\in\NN$}\quad\hbox{and}\quad
\widehat{h}^{{(j)}}_r(s)\ll r^{2\rho_j+1}.
\end{align*}
Recalling that $d_j=2\rho_j+1$,
let us define the following spectral estimator : 
$$
H_{r,R}({\bf s}(P^c), {\bf t}(P)):=\prod_{j=1}^k  r^{d_j}\prod_{j\in I_{pr}} (1+rt_j)^{-N} \prod_{j\in J_{pr}} (1+Rt_j)^{-N}  
\prod_{j\in J_{comp}} e^{2s_j(d_j-1)R} .
$$
Then it remains to analyze, for any given subset $P$ 
\begin{equation}\label{eq:sum}
\sum_{{\bf s}(P^c)\,,\,{\bf t}(P)} m\big((\otimes_{j\in P} \pi_{it_j})\otimes  (\otimes_{j\in P^c} \pi_{s_j}) ,\Gamma\big)
H_{r,R}({\bf s}(P^c),{\bf t}(P)).
\end{equation}

{\bf Part II : Multi-dimensional integration by parts}. It will be convenient to interpret this sum as a
multidimensional Riemann-Stieltjes integral, whose theory we now review (see \cite{fre}, \cite{z}). Let $\phi$ be a function defined on
$\prod_{j=1}^\ell [\alpha_j,\beta_j]$.
For $1\le i\le \ell$, we define
a function $\Delta_i\phi$ on 
the box $\prod_{j\ne i} [\alpha_j,\beta_j]$ by
$$
\Delta_i \phi=\phi|_{x_i=\beta_i}-\phi|_{x_i=\alpha_i}.
$$
Let 
$$
\Delta\phi=\left( {\prod}_{j=1}^\ell \Delta_j\right)\phi.
$$
Given a function $\phi$ defined on a box $B=\prod_{j=1}^\ell [a_j,b_j]$,
we consider partitions $B=\cup_{\cE}B({\cE})$ defined by arbitrary finite partitions
of the intervals $[a_j,b_j]$. The variation of $\phi$ 
in the sense of Vitali is defined as
$$
V(\phi)=\sup {\sum}_{\cE} \big|\Delta(\phi|_{B({\cE})})\big|,
$$
where the supremum is taken over all partitions as above.
Let us assume that 
$V(\phi)<\infty$ and $f$ is a continuous function
on $B=\prod_{j=1}^\ell [a_j,b_j]$.
Given a partition $B=\cup_{\cE} B({\cE})$, we choose points $\xi_{\cE}\in B({\cE})$.
The multidimentional Riemann-Stieltjes integral $\int_B f\, d\phi$ introduced by Fr\'echet in in \cite{fre} is defined as limit of the sums
$$
{\sum}_{\cE} f(\xi_{\cE}) \Delta(\phi|_{B({\cE})})
$$
as the diameter of the partitions goes to zero.
We note that when $\phi$ is continuously differentiable $\ell$ times, this integral can be expressed in terms of the usual Riemann integral:
$$
\int_B f\, d\phi =\int_B f\cdot \frac{\partial^\ell\phi}{\partial x_1\cdots\partial x_\ell}\, dx.
$$

We also need an integration by parts formula developed by Zaremba  \cite{z}. For this purpose,  we 
additionally assume that the restrictions of $\phi$
to the faces of $B$ also have bounded variation in the sense of Vitali. 
For a  subset $S$ of $\{1,\ldots,\ell\}$, we set
$B_S:=\prod_{j\in S} [a_j,b_j]$.
We define a function $\Delta_S\phi$ on  the box $B_{S^c}$ by
$$
\Delta_S \phi=\left({\prod}_{j\in S} \Delta_j\right)\phi.
$$
Using this notation, one has the following integration by parts
formula (see \cite[Prop.~2]{z}):
\begin{align}\label{eq:parts}
\int_B f(x)\, d\phi(x) &=\sum_S (-1)^{|S|}  \Delta_{S^c}\left( 
\int_{B_S} \phi(y)\, d_Sf(y)
\right) \\
&
=\sum_S (-1)^{|S|} 
\int_{B_S} \Delta_{S^c}\phi(y)\, d_Sf(y), \nonumber
\end{align}
where the sum is taken over all subsets $S$ of $\{1,\ldots,\ell\}$, and the notation $d_S$ designates that we take the Riemann-Stieltjes integral with respect to the variables indexed by $S$.

\medskip

We now return to analysis of the sums \eqref{eq:sum}, for any given subset $P\subset \set{1,\dots,\ell}$.
For ${\boldsymbol \sigma}=(\sigma_j)_{j\in P}$ and 
${\boldsymbol \tau}=(\tau_j)_{j\in P^c}$ with $\sigma_j\in [0,1/2)$ and $\tau_j\in [0,\infty)$,
let us denote by $N({\boldsymbol\sigma},{\boldsymbol\tau})$ the sum of multiplicities
$m\big((\otimes_{j\in P} \pi_{it_j})\otimes  (\otimes_{j\in P^c} \pi_{s_j}),\Gamma\big)$,
where ${\bf t}=(t_j)_{j\in P}$ runs over $\prod_{j\in P} (0,\tau_j]$
and ${\bf s}=(s_j)_{j\in P^c}$ runs over 
$\prod_{j\in P^c} [\sigma_j,1/2)$.
Then the above sum can be rewritten as a limit of  multidimensional Riemann-Stieltjes integrals:
\begin{equation}\label{eq:ii}
\pm \lim_{T\to\infty} \int_{[0,1/2]^{P^c}\times [0,T]^P} 
H_{r,R}({\boldsymbol\sigma},{\boldsymbol\tau})\, dN({\boldsymbol\sigma},{\boldsymbol\tau}).
\end{equation}
Applying \eqref{eq:parts}, the estimate reduces to
analysis of the expressions
\begin{align} \label{eq:int_int}
&\Delta_{S^c}\left(\int_{[0,1/2]^{S_{comp}}\times [0,T]^{S_{pr}}} N({\boldsymbol\sigma},{\boldsymbol\tau}) 
D_S(H_{r,R})({\boldsymbol\sigma},{\boldsymbol\tau})\, d_{S_{comp}}{\boldsymbol \sigma} d_{S_{pr}}{\boldsymbol \tau}\right)\\
=&
\int_{[0,1/2]^{S_{comp}}\times [0,T]^{S_{pr}}} \Delta_{S^c}\big(N({\boldsymbol\sigma},{\boldsymbol \tau}) 
D_S(H_{r,R})({\boldsymbol\sigma},{\boldsymbol \tau})\big)\, d_{S_{comp}}{\boldsymbol\sigma} d_{S_{pr}}{\boldsymbol\tau},\nonumber
\end{align}
where $S$ is a subset of $\{1,\ldots,\ell\}$,
$S_{pr}=S\cap P$, $S_{comp}=S\backslash P$,
and $D_S$
denotes the partial derivative of order $|S|$ with respect to the variables 
indexed by $S$. If $S\cap I_{comp}\ne \emptyset$, 
we get $D_S(H_{r,R})=0$, so that from now on we assume that $S\cap I_{comp}= \emptyset$. We obtain
\begin{align*}
D_S(H_{r,R})({\boldsymbol\sigma},{\boldsymbol \tau})
=\pm  \prod_{j=1}^k r^{d_j}  & \cdot
r^{|S\cap I_{pr}|} N^{|S\cap I_{pr}|} \prod_{j\in I_{pr}\cap S} (1+r\tau_j)^{-N-1} \prod_{j\in I_{pr}\backslash S} (1+r\tau_j)^{-N}\\
&\times 
R^{|S\cap J_{pr}|} N^{|S\cap J_{pr}|}\prod_{j\in J_{pr}\cap S} (1+R\tau_j)^{-N-1} \prod_{j\in J_{pr}\backslash S} (1+R\tau_j)^{-N} \\
&\times \prod_{j\in J_{comp}\cap S} (2(d_j-1)R) \prod_{j\in J_{comp}} e^{2\sigma_j(d_j-1)R}.
\end{align*}
We remark that $N({\boldsymbol\sigma},{\boldsymbol\tau})=0$ if one of the coordinates $\tau_j$ is equal to 0. Also 
$N({\boldsymbol\sigma},{\boldsymbol\tau})=0$ if one of the coordinates $\sigma_j$ is equal to $1/2$ because we have excluded the trivial representation.
Therefore,
\begin{align*}
\Delta_{S^c}\big(N({\boldsymbol\sigma},{\boldsymbol \tau}) D_S(H_{R,r})&({\boldsymbol\sigma},{\boldsymbol\tau})\big)\\
=&\pm N_S({\boldsymbol\sigma},{\boldsymbol\tau}; T)\\
&\times \prod_{j=1}^k r^{d_j}\cdot 
r^{|S\cap I_{pr}|} N^{|S\cap I_{pr}|} \prod_{j\in I_{pr}\cap S} (1+r\tau_j)^{-N-1} \prod_{j\in I_{pr}\backslash S} (1+rT)^{-N}\\
&\times 
R^{|S\cap J_{pr}|} N^{|S\cap J_{pr}|}\prod_{j\in J_{pr}\cap S} (1+R\tau_j)^{-N-1} \prod_{j\in J_{pr}\backslash S} (1+RT)^{-N} \\
&\times \prod_{j\in J_{comp}\cap S} (2(d_j-1)R) \prod_{j\in J_{comp}\cap S} e^{2\sigma_j(d_j-1)R}.
\end{align*}
where $N_S({\boldsymbol\sigma},{\boldsymbol\tau}; T)$ is obtained from 
$N({\boldsymbol\sigma},{\boldsymbol\tau})$ by substituting $\sigma_j=0$ for $j\in S^c\backslash P$ and $\tau_j=T$ for 
$j\in S^c\cap P$.

We also use the bound from Theorem \ref{density hypothesis}:
$$
N_S({\boldsymbol\sigma},{\boldsymbol\tau}; T) \ll_{\delta} \prod_{j\in P\cap S} (1+\tau_j)^{d_j(1-2m({\boldsymbol\sigma}))+\delta} \prod_{j\in P\backslash S} (1+T)^{d_j(1-2m({\boldsymbol\sigma}))+\delta}
$$
for all $\delta>0$, where $m({\boldsymbol\sigma}):=\max (\sigma_j:\, j\in S\backslash P)$. Now the computation reduces to one-dimensional integrals.
For $i\in I_{pr}\cap S$ and $r\in (0,1)$,
\begin{align*}
&\int_0^T (1+\tau_j)^{d_j(1-2m({\boldsymbol\sigma}))+\delta} (1+r\tau_j)^{-N-1} \, d\tau_j \\
\le &  r^{-1} \int_0^\infty  (1+\tau_j/r)^{d_j(1-2m({\boldsymbol\sigma}))+\delta} (1+\tau_j)^{-N-1} \, d\tau_j\\
\le &  r^{d_j(2m({\boldsymbol\sigma})-1)-1-\delta} \int_0^\infty  (1+\tau_j)^{-N+2} \, d\tau_j \ll_N  r^{d_j(2m({\boldsymbol\sigma})-1)-1-\delta}
\end{align*}
when $N\ge  4$.
Also for $j\in J_{pr}\cap S$ and $R> 1$,
$$
\int_0^T (1+\tau_j)^{d_j(1-2m({\boldsymbol\sigma}))+\delta} (1+R\tau_j)^{-N-1} \, d\tau_j \le \int_0^\infty  (1+\tau_j)^{-N+2} \, d\tau_j\ll_N 1
$$
when $N\ge  4$.
It follows that ultimately, up to a constant (depending on $N$ and $\delta$), we get the bound:
\begin{align*}
&\prod_{j=1}^k r^{d_j}\cdot 
r^{|S\cap I_{pr}|}  \cdot  \prod_{j\in I_{pr}\cap S} r^{d_j(2m({\boldsymbol\sigma})-1)-1-\delta}\cdot  (1+rT)^{-N|I_{pr}\backslash S|}\\
&\times
R^{|S\cap J_{pr}|} (1+R T)^{-N|J_{pr}\backslash S|} \cdot R^{|S\cap J_{comp}|}\prod_{j\in J_{comp}\cap S} e^{2m({\boldsymbol\sigma})(d_j-1)R} 
\prod_{j\in P\backslash S} (1+T)^{d_j(1-2m({\boldsymbol\sigma}))+\delta}.
\end{align*}
If either $I_{pr}\backslash S\ne \emptyset$ or 
$J_{pr}\backslash S\ne \emptyset$, this clearly vanishes as $T\to \infty$,
provided that $N$ is sufficiently large.
Hence, we may assume that $I_{pr}\subset S$ and $J_{pr}\subset S$. Then also $P\backslash S=\emptyset$,
and we simplify this bound as 
\begin{align*}
& \prod_{j=1}^k r^{d_j}\cdot 
  \prod_{j\in I_{pr}\cap S} r^{d_j(2m({\boldsymbol\sigma})-1)-\delta}\cdot  R^{|S\cap J_{pr}|+|S\cap J_{comp}|}\cdot \prod_{j\in J_{comp}\cap S} e^{2m({\boldsymbol\sigma})(d_j-1)R} \\
\le & \prod_{j=1}^k r^{d_j}\cdot 
  \prod_{j=1}^k r^{d_j(2m({\boldsymbol\sigma})-1)-\delta}\cdot  R^{\ell-k}\cdot \prod_{j=k+1}^\ell e^{2m({\boldsymbol\sigma})(d_j-1)R} \\
 \ll &_{\delta}\,\, R^{\ell-k} r^{-k \delta}  \left(\prod_{j=1}^k r^{d_j}\right)^{2m({\boldsymbol\sigma})} 
 \left(\prod_{j=k+1}^\ell e^{(d_j-1)R}\right)^{2m({\boldsymbol\sigma})} 
\end{align*}
for all $\delta>0$,
where we used that $I_{pr}\subset \{1,\ldots,k\}$ 
and $J_{pr}\cup J_{comp}= \{k+1,\ldots,\ell\}$. 
This provides an estimate for the part of 
the integral \eqref{eq:int_int} where have integrated over $\boldsymbol\tau$.
Now it remains to integrate over ${\boldsymbol\sigma}\in [0,1/2]^{S_{comp}}$.
Taking the (weaker form of the) matching volumes condition stated as \eqref{eq:rr} into account, we deduce that the last sum is bounded by $\ll_{\delta} e^{c\delta R}$
with fixed explicit $c>0$, which is the 
required estimate proving Theorem  \ref{trace-estimate}. 
\end{proof}

\begin{remark}
 Interestingly, the worst case arises only when
$$
I_{pr}=\{1,\ldots,k\}\quad\hbox{and}\quad  J_{comp}=\{k+1,\ldots,\ell\}
$$
and the other cases all satisfy the stronger bound of the form $O(e^{-\alpha R})$
with some $\alpha>0$.
\end{remark}

\subsection{Completion of the proof of Theorem \ref{th:quaternion}}
\label{sec:proof}

We deduce Theorem \ref{th:quaternion} from Theorem \ref{th:general}.   Let 
$$
G_1:=\prod_{j=1}^k G^{(j)}\quad\hbox{and}\quad 
G_2:=\prod_{j=k+1}^\ell G^{(j)},
$$
and
$$
K_1:=\prod_{j=1}^k K^{(j)}\quad\hbox{and}\quad 
K_2:=\prod_{j=k+1}^\ell K^{(j)},
$$
where $K^{(j)}$ is a maximal compact subgroup of $G^{(j)}$.
For $k+1\le j\le \ell$, we choose $x_0^{(j)}$ so that $K^{(j)}=\hbox{Stab}_{G^{(j)}}(x_0^{(j)})$.
We recall the metrics on $G_1/K_1$ and $G_2/K_2$ defined by letting $D^{(j)}$ denote the Riemannian metric of constant curvature $-1$ on $G^{(j)}/K^{j)}$, and 
\begin{align*}
{D}_1(x_1,y_1)&:=\max \Big(D^{(j)}(x_1^{(j)},y_1^{(j)}):\, j=1,\ldots, k\Big)\quad  \hbox{ for  
$x_1,y_1\in G_1/K_1,$}\\
{D}_2(x_2,y_2)&:=\max \Big(D^{(j)}(x_2^{(j)},y_2^{(j)}):\, j=k+1,\ldots, \ell\Big)\quad \hbox{ for $x_2,y_2\in G_2/K_2.$}
\end{align*}
The sets $\widetilde{B}^{G_1}_{r}\subset G_1$ and
$\widetilde{B}^{G_2}_{R}\subset G_2$ are defined as in \eqref{eq:ball}.
Note that under the foregoing definitions, the volume growth exponent of the balls $\widetilde{B}^{G_2}_{R}\subset G_2$ is indeed $a_{[k+1,\ell]}$.

Now we proceed to verify the conditions of Theorem \ref{th:general}. 
Let $r\in (0,1)$, $R\in (1,\infty)$, $\delta \in (0,1)$ and $\delta^\prime\in (0,1)$.
To verify condition (T), we use functions
$f_{1,r}^{(\rho_j)}\in C_c(G^{(j)})$ with $j=1,\ldots, k$
and $f_{2,R}^{(\rho_j)}\in C_c(G^{(j)})$ with $j=k+1,\ldots, \ell$ constructed in Lemmas \ref{l:small} and \ref{l:large}. We recall that these families of functions depend upon a bump function $\eta$, supported on $(-1+\delta^\prime, 1-\delta^\prime)$. 

We define functions
$$
f_{1,r}:=f_{1,r}^{(\rho_1)}\otimes \cdots \otimes f_{1,r}^{(\rho_k)}\in C_c(G_1)\quad\hbox{and}\quad
f_{2,R}:=f_{2,R}^{(\rho_{k+1})}\otimes \cdots \otimes f_{2,R}^{(\rho_\ell)}\in C_c(G_2)
$$
and
$$
F_{r,R}:=f_{1,r}\otimes  f_{2,R}\in C_c(G_1\times G_2).
$$
As above let $d_j=2\rho_j+1$. We recall that these function satisfy:
\begin{itemize}
\item  For $j=1,\ldots, k$  and $r\in (0,1)$,
$$
\supp(f^{(\rho_j)}_{1,r})\subset \widetilde B_{r}^{G^{(j)}}  \quad \hbox{and} \quad  \int_{G^{(j)}} f^{(\rho_j)}_{1,r}\, dm_{G^{(j)}} \asymp_{\,\,\delta^\prime} r^{d_j/2}
$$
and 
for $z=\sigma+i\tau \in \widehat{G^{(j)}}^{sph}$,
\begin{align*}
\abs{\widehat{f^{(\rho_j)}_{1,r}}(z)}\ll_{N} r^{d_j/2}(1+r\abs{\tau})^{-N}\;\hbox{ for all $N\in\NN$}.
\end{align*}

\item For $j=k+1,\ldots, \ell$ and $R\in (1,\infty)$,
$$
\supp(f^{(\rho_j)}_{2,R})\subset \widetilde B_R^{G^{(j)}}\quad\hbox{and}\quad e^{(1-\delta^\prime)(d_j-1)R/2}\ll_{\,\,\delta^\prime} \int_{G^{(j)}} f^{(\rho_j)}_{2,R}\, dm_{G^{(j)}} \ll e^{(d_j-1)R/2}
$$
and for $z=\sigma+i\tau \in \widehat{G^{(j)}}^{sph}$,
\begin{align*}
\abs{\widehat{f_{2,R}^{(\rho_j)}}(z)}\ll_{N} e^{\sigma (d_j-1)R}(1+R\abs{\tau})^{-N}\;\hbox{ for all $N\in\NN$}.
\end{align*}
\end{itemize}

Therefore, 
\begin{align}
\supp( f_{1,r})&\subset  \widetilde{B}^{G_1}_{r}\quad\hbox{and}\quad \supp(f_{2,R})\subset \widetilde{B}^{G_2}_{R}, \nonumber\\
\int_{G_1} {f_{1,r}}\, dm_{G_1}&=\prod_{j=1}^k \int_{G^{(j)}} f_{1,r}^{(\rho_j)}\, dm_{G^{(j)}}\asymp_{\,\, \delta^\prime} \prod_{j=1}^k r^{d_j/2}, \nonumber\\
\prod_{j=1}^k e^{(1-\delta^\prime)(d_j-1)R/2}&\ll_{\delta^\prime} \int_{G_2} {f_{2,R}}\, dm_{G_2}=\prod_{j=k+1}^\ell \int_{G^{(j)}} f_R^{(j)}\, dm_{G^{(j)}}\ll  \prod_{j=1}^k e^{(d_j-1)R/2}.
\label{eq:low_vol}
\end{align}
Let us define $r={\bf r}(R)$ by the condition

\begin{equation}\label{eq:v}
\prod_{j=1}^k r^{d_j/2}=\prod_{j=1}^k e^{-(d_j-1)R/2}.
\end{equation}
Namely,
$$
{\bf r}(R):= \exp\left({-\frac{a_{[k+1,\ell]}}{d_{[1,k]}}R}\right).
$$
We note that in view of volume estimate \eqref{volume-SL2},
this choice of $r$ correspond to the matching volume condition:
$$
m_{G_1}(\widetilde{B}_r^{G_1})\cdot m_{G_2}(\widetilde{B}_R^{G_2})\asymp 1.
$$
Hence, we can apply Theorem \ref{trace-estimate} to 
conclude that  
$$
\hbox{\rm Tr}\Big(\lambda_{\Gamma\backslash G}\left({F}^\ast_{{\bf r}(R),R}\ast{F}_{{\bf r}(R),R}\right) \Big) \ll_\delta e^{c\delta R}
$$
for some fixed $c> 0$.
Additionally, we rescale the functions 
$f_{2,R}^{(\rho_j)}$, $j=k+1,\ldots,\ell$
to arrange that
$$
\int_{G_1\times G_2} {F}_{{\bf r}(R),R}\, dm_G =\left(\int_{G_1} {f}_{1,{\bf r}(R)}\, dm_{G_1}\right)\left(\int_{G_2} {f}_{2,R}\, dm_{G_2}\right)=1.
$$
In view of the estimate \eqref{eq:low_vol},
this rescaled function still has the trace bound of the form
$$
\hbox{\rm Tr}\Big(\lambda_{\Gamma\backslash G}\left({F}^\ast_{{\bf r}(R),R}\ast{F}_{{\bf r}(R),R}\right) \Big) \ll_{\delta, \delta^\prime} e^{(c\delta+c^\prime\delta^\prime) R}
$$
for some fixed $c^\prime> 0$ independent of $\delta$, $\delta^\prime$ and $R$. Given $\delta > 0$, we can choose say $\delta^\prime=\delta$, and 
this verifies condition (T).

To verify the condition (SP), 
we consider the unitary representation
$\sigma$ of $G_2$ on $L^2_0(\Gamma\backslash G)$,
the orthogonal complement of constant functions.
It follows from Theorem \ref{strong SP} and \eqref{integrability} that this representation has positive
integrability exponent. Hence, it follows from \cite{n}
that some even tensor power $\sigma^{\otimes n}$ is weakly
contained in the regular representation $\lambda_{G_2}$ on 
$L^2(G_2)$, and for all non-negative $\psi\in L^1(G_2)$,
$$
\norm{\sigma(\psi)}\le \norm{\lambda_{G_2}(\psi)}^{1/n}.
$$
Furthermore, by the Kunze--Stein inequality \cite{c78},
$$
\norm{\lambda_{G_2}(\psi)}\ll_p \norm{\psi}_p\quad
\hbox{for $p\in [1,2)$.}
$$
Fix $c_0>0$ and consider the function
$$
\chi_R:= m_{G_2}(\widetilde{B}^{G_2}_{R-c_0})^{-1}\cdot \chi_{\widetilde{B}^{G_2}_{R-c_0}}.
$$
Let $\phi\in C_c(G_2)$ be a $K_2$-bi-invariant function
such that 
$$
\supp(\phi)\subset \widetilde{B}^{G_2}_{c_0}\quad
\hbox{and}\quad
\int_{G_2}\phi\, dm_{G_2}=1
$$
and set 
$$
\psi_R:=\phi*\chi_R.
$$
This is a continuous $K_2$-bi-invariant non-negative function such that
$$
\supp(\psi_R)\subset \widetilde{B}^{G_2}_{R}\quad\hbox{and}\quad
\int_{G_2}\psi_R\, dm_{G_2}=1.
$$
Also, fixing $p$ satusfying $1< p < 2$ :
$$
\norm{\sigma(\psi_R)}\ll \norm{\sigma(\chi_R)}\ll \norm{\chi_R}_p^{1/n}\ll m_{G_2}(\widetilde{B}^{G_2}_{R-c_0})^{-(1-1/p)/n}\ll e^{-\mathfrak{w} R}
$$
for some fixed $\mathfrak{w} >0$, with the impied constant independent of $R$. This verifies the condition (SG).

Now we can apply Theorem \ref{th:general} to conclude that 
for every $\delta>0$ and almost every $(g_1,g^\prime_1)\in G_1\times G_1$,
the inequalities 
$$
D_1(\gamma_1 g_1K_1,g^\prime_1K_1)\le  2{\bf r}(R)\quad\hbox{and}\quad D_2(\gamma_2 K_2,K_2) \le (1+\delta)R
$$
have solutions $\gamma=(\gamma_1,\gamma_2)\in\Gamma$ for all $R>R_0(g_1,g^\prime_1,\delta)$.
This implies that for every $\zeta> \frac{d_{[1,k]}}{a_{[k+1,\ell]}}$ and almost every $x=(x^{(1)},\ldots,x^{(k)}), y=(y^{(1)},\ldots,y^{(k)})\in {X}_{[1,k]}$, the inequalities
\begin{align*}
&D^{(j)}(\gamma^{(j)}\cdot x^{(j)},y^{(j)})\le \epsilon ,\;\; j=1,\ldots,k,\\
&D^{(j)}(\gamma^{(j)}\cdot x_0^{(j)},x_0^{(j)})\le \zeta \, \log(1/\epsilon), \quad j=k+1,\ldots,\ell,
\end{align*}
have solutions $\gamma=(\gamma^{(1)},\ldots,\gamma^{(\ell)})\in\Gamma$ for all $\epsilon\in \big(0,  \epsilon_0(x,y,\zeta)\big)$.  

Note that for $j=1,\dots,k$ the first $k$ inequalities for $\epsilon \le 1$ (say), imply that the projection of $\gamma$ to $G^{(j)}$ is contained in a fixed bounded set. Therefore, for almost every $x,y\in {X}_{[1,k]}$,
and for the metric $\cD$ defined on $G$ (see \ref{metricD})), the solutions $\gamma$ to the system of inequalities satisfy $\cD(e,\gamma)\le \zeta \log (1/\epsilon)+C$. 
It follows that for the metric $\cD$ on $G$ we have the following bound on the Diophantine approximation exponent : 
\begin{equation}\label{eq:llll}
\kappa_\Gamma(x,y)\le \frac{d_{[1,k]}}{a_{[k+1,\ell]}}.
\end{equation}
Furthermore, we note that 
$\dim ({X}_{[1,k]})={d_{[1,k]}}$, and
it follows from \eqref{volume-SL2} that
$$
m_{G_2}\big(\widetilde{B}_R^{G_2}\big)\asymp \exp\big({a_{[k+1,\ell]} R}\big)
\quad\hbox{as $R\to\infty$.}
$$
Thus, it follows from a pigeonhole argument
as in \cite[Th.~3.1]{GGN14} that, in fact, the equality in \eqref{eq:llll} holds 
for almost every $x,y\in {X}_{[1,k]}$.

This concludes the proof of Theorem  \ref{th:quaternion}.\qed

\section{Local estimates for general spherical transforms}\label{sec:local SF}

Let $G$ be a connected semisimple real Lie group with finite center.
We fix a Cartan involution $\Theta$ of $G$ and denote by $K$ the 
corresponding maximal compact subgroup. 
Let $\mathfrak s$ be the $-1$-eigenspace of in the Lie algebra of $G$ 
corresponding to $\Theta$. The map
$$
K\times \mathfrak s\to G: (k,Z)\mapsto k\exp(Z)
$$
is a diffeomorphism. 
Let $w \colon \mathfrak s\to \mathbb R_{\geq 0}$ be a smooth compactly supported function
which is invariant under $\hbox{Ad}(K)$.
We define $\omega: G:\to \mathbb R_{\geq 0}$ by
$$
\omega(ke^Z):=w\left(Z\right)\quad\hbox{for $(k,Z)\in K\times \mathfrak s$}.
$$
This is a smooth compactly supported $K$-bi-invariant function on $G$.
Our goal is to estimate the spherical transform of the function $\omega$.
In fact, we will be interested in the following family
of rescalings of $\omega$ defined  for $r\in (0,1]$ by
$$
\omega_r(ke^Z):=r^{-d} w\left(r^{-1} Z\right)\quad\hbox{for $(k,Z)\in K\times \mathfrak s$},
$$
where $d=\dim(\mathfrak s)=\dim (G)-\dim (K)$.
Equivalently, these functions can be defined in terms of a Cartan subalgebra
$\mathfrak a \subset \mathfrak s$. Then
$$
G=K\exp(\mathfrak a)K,
$$
and 
\begin{equation}\label{omega}
\omega_r(k_1e^H k_2)=\omega_r\left(k_1k^{-1}_2e^{\text{Ad} \,k_2(H)}\right)=
r^{-d} w\left(r^{-1}H\right)\quad\hbox{for $k_1,k_2\in K$ and $H\in \mathfrak a$}.
\end{equation}

We fix a minimal parabolic subgroup $P$ of $G$ containing $A$, and let $N$ be the unipotent radical of $P$. The the product map $K\times \exp(\mathfrak a)\times N\to G$ is a diffeomorphism.
For $g\in G$, we write $g=k_ge^{H_I(g)}n_g$ with $k_g\in K$, $H_I(g)\in \mathfrak a$, and 
$n_g\in N$. 

{\it Note about notation:} In this section $\sigma, \tau, \rho$ will denote vectors in $\text{\rm Hom}_\RR(\mathfrak{a},\RR)$, generalizing the notation used in previous sections when $A$ was assumed one-dimensional. 

We recall that for any $\xi\in \mathfrak a^*_\mathbb C:=\text{Hom}_\RR(\mathfrak{a},\CC)$, the elementary spherical function $\varphi_\xi$ is defined by 
$$\varphi_\xi(g):=\int_K e^{(\xi-\rho)(H_I(gk))}\, dm_K(k),$$ where $m_K$ is the probability Haar measure on $K$, and $\rho$ denotes half the sum of the positive roots. 
The spherical transform of $\omega_r$ is defined by
$$
\widehat{\omega}_r(\xi):= \int_G \omega_r(g)\varphi_{-\xi}(g)\,dm_G(g) \quad \hbox{for $\xi\in \mathfrak a^*_{\mathbb C}.$}
$$

\begin{theorem}\label{lem-1} For any $N\geq 1$ and $B>0$ there exists a constant $C_{N,B}>0$ such that for  all $\sigma,\tau\in \mathfrak a^*:=\text{\rm Hom}_\RR(\mathfrak{a},\RR)$ with $ \|\sigma\|\leq B$
and  $r\in (0,1]$,
\begin{equation}\label{local-spherical}  \abs{\widehat{\omega}_r(\sigma+i\tau)}\leq C_{N,B} (1+r\|\tau\|)^{-N}.
\end{equation} 
\end{theorem}

First, we observe that since $\omega_r$ is right-$K$-invariant,
\begin{align}
 \widehat{\omega}_r(\xi) &=
 \int_G \omega_r(g)\int_K  e^{-(\xi+\rho)(H_I( gk))}\, dm_K(k) dm_G(g) \nonumber \\
&=\int_G\int_K \omega_r(gk )  e^{-(\xi+\rho)(H_I( gk))}\, dm_K(k) dm_G(g) \nonumber\\
&=\int_K\left(\int_G \omega_r(gk )  e^{-(\xi+\rho)(H_I( gk))}\,  dm_G(g)\right) dm_K(k)\nonumber\\
&= 
\int_G \omega_r(g)e^{-(\xi+\rho)(H_I(g))}\,dm_G(g). \label{eq:iii}
\end{align}
Let us now use the Iwasawa decomposition of $G$ in a different form, namely  $G=KN\exp(\mathfrak a)$.
Haar measure on $G$ is given in these coordinates by
\begin{equation} \label{Haar} 
\int_G f(g)\,dm_G(g)=\int_{K\times{\mathfrak n}\times {\mathfrak a}} f(ke^Y e^H)\, dm_K(k) dY dH, \quad f\in C_c(G),
\end{equation} 
where $\mathfrak n$ denotes the Lie algebra of $N$, and $dY$ and  $dH$
are the Lebesgue measure on $\mathfrak n$ and on $\mathfrak a$ respectively.
For $g\in G$, we write  $g=k_g^\prime n_g^\prime e^{H_I^\prime(g)}$ with 
$k_g^\prime\in K$,  $n_g^\prime\in N$,  $H_I^\prime(g)\in \mathfrak a$. The coordinates here are of course different from those  defining the component $H_I(g)$ appearing above, but 
it follows from the uniqueness of the Iwasawa decomposition that $k_g'=k_g$, $H_I^\prime(g)=H_I(g)$, and 
$n_g'=e^{H_I(g)} n_g e^{-H_I(g)}$.
Combining (\ref{Haar}) with (\ref{eq:iii}) and using the left-$K$-invariance of $\omega_r,$
we conclude that
\begin{equation}\label{kernel} \widehat{\omega}_r(\xi)=\int_{{\mathfrak n}\times {\mathfrak a}} \omega_r(e^{Y}e^{H})e^{-(\xi+\rho)(H)}\,dYdH.
\end{equation} 

The foregoing expression is the Euclidean Fourier-Laplace transform of a family of local bump functions  on the Euclidean space ${\mathfrak n}\times {\mathfrak a}$, and we  proceed to analyze its behavior.  
Let us introduce the diffeomorphism
$$
\Psi\colon \mathfrak n \times \mathfrak a\to \mathfrak s: (Y,H)\mapsto \cJ^{-1} (e^Ye^HK),
$$
where $\cJ$ is the diffeomorphism 
$$
\cJ\colon \mathfrak s\to G/K: Z\mapsto \exp(Z)K.
$$
We note that $\Psi(0)=0$.  Then 
\begin{align*}
\widehat{\omega}_r(r^{-1}\xi)&=r^{-d}\int_{{\mathfrak n}\times {\mathfrak a}}w\left(r^{-1}\Psi(Y,H)\right)e^{-r^{-1}(\xi+r\rho)(H)}\,dYdH\\    
&=\int_{{\mathfrak n}\times {\mathfrak a}}w\left(\Psi_r(Y,H)\right)e^{-(\xi+r\rho)(H)}\,dYdH\,,
\end{align*}
where the maps $\Psi_r$ are defined by
$$
\Psi_r\colon \mathfrak n\times \mathfrak a\to \mathfrak s:
(Y,H)\mapsto  r^{-1}\Psi(rY,rH).
$$
Here we have used the change of variables 
$(Y,H)\mapsto (rY', rH')$, noting that $d =\dim({\mathfrak a})+\dim ({\mathfrak n })$.  Letting $\mathcal L$ denote the Euclidean Fourier-Laplace transform on $\mathfrak n\times \mathfrak a\simeq \mathbb R^d$, we have 
\begin{equation}\label{composition} \widehat{\omega}_r(r^{-1}\xi)=\mathcal L(w\circ \Psi_r)\big(-\xi-r\rho\big).
\end{equation} 
 
 The proof of Theorem \ref{lem-1} will utilize the following Lemma which we now turn to prove. 
\begin{lemma}\label{lem-2}
Let $w\colon \mathbb R^d\to \mathbb R$ be a smooth compactly supported function and let $\Psi\colon \mathbb R^d\to \mathbb R^d$ be a diffeomorphism such that $\Psi(0)=0$. Write $\Psi_r\colon \mathbb R^d\to\mathbb R^d$ for the map $v\mapsto r^{-1}\Psi(rv)$ when $r \neq 0$, and $\Psi_0(v)=D\Psi(0)(v)$.  Then, there exists a compact set $\cW\subset \mathbb R^d$ such that $w\circ \Psi_r$ is supported in $\cW$ for any $r\in [-1,1]$ and for every $N$ there exists a constant $C_N$ such that the $N$-th partial derivatives of $w\circ \Psi_r$ are bounded by $C_N$, independently of $r\in [-1,1].$
\end{lemma}

\begin{proof}

For any $v\in \mathbb R^d$ let $Q(v)\colon \mathbb R^d\times \mathbb R^d\to\mathbb R^d$ be the $\RR^d$-valued bilinear map defined by 
$$Q(v)(u_1,u_2):=\int_0^1 (1-s) D^2\Psi(sv)(u_1,u_2)ds.$$
Here  $D^2\Psi$ is the second derivative of $\Psi$ so that at any point it takes values in bilinear forms $\RR^d\times \RR^d\to \RR^d$.
Clearly $Q$ is a smooth function of $v\in \RR^d$. Consider the function $y(t) =\Psi(tv)$, so that $y'(t)=D\Psi(tv)(v)$ and $y''(t)=D^2\Psi(tv)(v)(v)$. Applying Taylor's expansion to first order with integral form of the remainder to it at $t=1$, the expression $y(1)=y(0)+y^\prime(0) +\int_0^1(1-s)y''(s)ds$ translates to 
$$\Psi(v)=\Psi(0)+D\Psi(0)(v)+ Q(v)(v,v)=D\Psi(0)(v)+Q(v)(v,v).$$ 

Recalling that we defined $ \Psi_0(v)=D\Psi(0)(v)$,  we have for all $r\in \RR$
$$\Psi_r(v)=D\Psi(0)(v)+rQ(rv)(v,v)=r^{-1}\Psi(rv) \,\,(\text{when } r\neq 0).$$ 
Let us consider the smooth map
$$
\tilde \Psi: \RR^d \times \RR \to \RR^d: (v,r)\mapsto D\Psi(0)(v)+rQ(rv)(v,v).
$$

\begin{claim}\label{compact-support} 
 The map $\tilde\Psi |_{\mathbb R^d\times [-1,1]}$ is proper.
 \end{claim}

\begin{proof} Let $\cK\subset \RR^d$ be a compact set, and without loss of generality we may assume that $\cK$ is a convex symmetric set containing a ball centered at $0$, and show that $\tilde{\Psi}^{-1}(\cK)\cap\left(\RR^d\times [-1,1]\right)$ is bounded. Since $\Psi$ is a diffeomorphism fixing $0$, fix $ 0< c < 1$ such that $\Psi^{-1}  : c\cK\to \Psi^{-1}(c\cK)$  is Lipschitz with constant $C>0$. 
Consider $(v,r)\in \RR^d\times \RR $ such that $-1\le r \le 1$ and in addition $\tilde{\Psi}(v,r) \in \cK$, so that $(v,r)\in \tilde{\Psi}^{-1}(\cK)\cap\left(\RR^d\times [-1,1]\right)$. (For $r=0$ the expression is $D\Psi(0)(v)$, by the foregoing formula for $\tilde{\Psi}$).  It suffices to consider $r\ge 0$, and then either $r \in [c,1]$ or $r\in [0, c]$.   In the first case, since $\Psi(rv)\in r\cK\subseteq \cK$ we have  $v\in r^{-1}\Psi^{-1}(\cK)\subseteq 
[1,c^{-1}]\cdot \Psi^{-1}(\cK)$,  so that 
$$\tilde\Psi^{-1}(\cK)\cap \left(\mathbb R^d\times [c,1]\right) \subseteq  [1,c^{-1}]\cdot \Psi^{-1}(\cK)\times [c,1].$$ 
In the second case, when $ r > 0$, $\tilde{\Psi}(r,v)=r^{-1} \Psi(rv)\in \cK $ so that $\Psi(rv) \in r\cK\subseteq  c\cK$, $\Psi^{-1}$ has Lipschitz constant $C$ on $c\cK$, and so 
$$\norm{v}=r^{-1}\norm{rv}=r^{-1}\norm{\Psi^{-1}(\Psi(rv))}\le C \norm{r^{-1}\Psi(rv)}
\le C\max_{u\in \cK}\norm{u}:=C_{\cK}$$
 Therefore if $\cB_{\cK}$ denotes a closed ball centered at $0$ of radius $C_{\cK}$ then 
$$\tilde\Psi^{-1}(\cK)\cap \left(\mathbb R^d\times (0,c]\right)\subset \cB_{\cK}\times (0,c]\,.
 $$
Finally, when $r=0$, the linear map $\tilde\Psi(v,0)=\Psi_0(v)=D\Psi(0)(v)$ is proper, since $D\Psi(0)$ is a linear isomorphism, being the derivative of a diffeomorphism. 

Altogether $\tilde\Psi^{-1}(\cK)\cap \left(\mathbb R^d\times [-1,1]\right)$ is contained in a finite union of  compact sets so it is bounded and the claim is proved. 
\end{proof}

Now note that it follows from the claim that $w\circ \tilde \Psi$ restricted to $\mathbb R^d\times [-1,1]$ is compactly supported. Let $\cW$ be the projection of the support of $w\circ \tilde \Psi$ restricted to $\mathbb R^d\times [-1,1]$ onto $\mathbb R^d$. Then $\cW$ is compact and contains the support of $w\circ \Psi_r$ for every $r\in [-1,1].$ Finally, the assertion on the  derivatives follows from the fact that  $w\circ \tilde\Psi$ is smooth so the $N$-th partial derivatives are uniformly bounded on $\cW\times [-1,1].$ 
\end{proof}

\begin{proof}[Proof of  Theorem \ref{lem-1}] Fix now $\sigma+i\tau=\eta$ where $\sigma,\tau\in \mathfrak a^*:=\text{\rm Hom}_\RR(A,\RR)$, and also $ r \in (0,1]$. To estimate $\widehat{\omega}_r(\sigma+i\tau)$ write $\eta=r^{-1}\xi$. Then we use (\ref{composition})
to write 
$$\widehat{\omega}_r(\sigma+i\tau)=\widehat{\omega}_r(r^{-1}\xi)=
\mathcal L(w\circ \Psi_r)\big(-\xi-r\rho\big)=
\mathcal L(w\circ \Psi_r)\big(-r(\eta+\rho)\big).
$$
By Lemma \ref{lem-2}, $w\circ \Psi_r$ is smooth and compactly supported 
with uniformly bounded partial derivatives, on the set $\cW\times [-1,1]$. 
If $\norm{\sigma}\le B$, $r\in (0,1]$ and $r\norm{\tau} \le 4(1+\norm{\rho}+B)$ (say), 
then $\widehat{\omega}_r(\sigma+i\tau)\le C_{1,B}$ which is all that Lemma \ref{lem-1} asserts in this case. Otherwise, let $\norm{\sigma}\le B$, $r \in (0,1]$ and $r\norm{\tau}\ge 4(1+\norm{\rho}+B)$, and in that case 
\begin{align*}
r\norm{\eta+\rho}&=r\norm{(\sigma+\rho)+i\tau} \ge \frac{1}{2} \norm{r\tau}+( \frac12 \norm{r\tau}-r\norm{\sigma+\rho} )\\
&\ge \frac12\norm{r\tau}+1\ge \frac12(1+r\norm{\tau}).
\end{align*}
Integrating the Laplace transform by parts once, we conclude that
$$\abs{\mathcal L(w\circ \Psi_r)(-r(\eta+\rho))}\le C'_{1,B}\norm{r(\eta+\rho)}^{-1}\le 
C'_{1,B}\left(1+r\norm{\tau}\right)^{-1}.
$$
The same argument is valid for repeated integrations by part, and this concludes the proof of Theorem \ref{lem-1}.
\end{proof}

\vskip2truein

Mikołaj Frączyk

Faculty of Mathematics and Computer Science, Jagiellonian University

mikolaj.fraczyk@uj.edu.pl

\vskip0.3truein

Alexanger Gorodnik 

Institute f\"ur Mathematik, Universit\"at Z\"urich 

alexander.gorodnik@math.uzh.ch

\vskip0.3truein 

Amos Nevo

Department of Mathematics, Technion, and 

Department of Mathematics, University of Chicago

Amosnevo6@gmail.com
\end{document}